\documentclass[reqno,makeidx,11pt]{amsart}
\usepackage{amssymb,amsfonts,amsmath,graphicx
}

\def\C{{\mathbb{C}}}

\def\E{{\mathbb{E}}}

\def\N{{\mathbb{N}}}
\def\0{{\mathbb{O}}}

\def\R{{\mathbb{R}}}
\def\T{{\mathbb{T}}}

\def\V{{\mathbb{V}}}
\def\Z{{\mathbb{Z}}}

\def\cB{{\mathcal B}}

\def\cH{{\mathcal H}}
\def\cI{{\mathcal I}}
\def\cK{{\mathcal K}}
\def\cL{{\mathcal L}}

\def\cP{{\mathcal P}}
\def\cR{{\mathcal R}}

\def\Id{{\rm Id \, }}

\def\sp{{\rm span \, }}
\def\Tr{{\rm Tr  }}

\newcommand{\norm}[1]{{\left\|{#1}\right\|}}

\newcommand{\abs}[1]{{\left|{#1}\right|}}
\newcommand{\scal}[1]{{\left\langle{#1}\right\rangle}}
\newcommand{\set}[1]{{\left\{{#1}\right\}}}
\newcommand{\un}{{\mathbf 1}}

\newcommand{\actson}{\curvearrowright}
\newcommand{\rd}{{\,\mathrm d}}
\newcommand{\wh}{\widehat}
\newcommand{\ovtimes}{\overline{\otimes}}

\newtheorem{thm}{Theorem}[section]
\newtheorem{cor}[thm]{Corollary}
\newtheorem{lem}[thm]{Lemma}
\newtheorem{prop}[thm]{Proposition}

\newtheorem*{thm*}{Theorem}
\newtheorem*{lem*}{Lemma}
\newtheorem*{cor*}{Corollary}

\theoremstyle{definition}
\newtheorem{defn}[thm]{Definition}

\newtheorem{exs}[thm]{Examples}

\theoremstyle{remark}
\newtheorem{rem}[thm]{Remark}

\title[Treeability and the Haagerup property for  groupoids]
{Old and new about  treeability and the Haagerup property for measured groupoids}

\author{Claire Anantharaman-Delaroche}
\address{Laboratoire de Math\'ematiques-Analyse, Probabilit\'es, Mod\'elisation-Orl\'eans (MAPMO - UMR6628),Ê
F\'ed\'eration Denis Poisson (FDP - FR2964),
CNRS/Universit\'e d'Orl\'eans,
B. P. 6759, F-45067 Orl\'eans Cedex 2}
\email{claire.anantharaman@univ-orleans.fr}



\begin{document}

\begin{abstract} This is mainly an expository text on the Haagerup property for countable groupoids equipped with a quasi-invariant measure,
aiming to complete an article of Jolissaint devoted to the study of this property for probability measure preserving countable equivalence relations.
We show that our definition is equivalent to the one given by Ueda in terms of the associated inclusion of von Neumann algebras.
It makes obvious the fact that treeability implies the Haagerup property for such groupoids. For the sake of completeness, we also describe, or recall, the connections
with amenability and Kazhdan property (T).
\end{abstract}
\maketitle

\section*{Introduction} 
Since the seminal paper of Haagerup \cite{Haa78}, showing that free groups have the (now so-called) Haagerup property, or property (H),
this notion plays an increasingly important role in group theory (see the book \cite{CCJJV}). A similar property (H) has been
introduced for finite von Neumann algebras \cite{Connes82, Choda} and it was proved in \cite{Choda}  that a countable group $\Gamma$
has property (H) if and only if its von Neumann algebra $L(\Gamma)$ has property (H).

 Later, given a von Neumann subalgebra $A$
of a finite von Neumann algebra $M$, a property (H) for $M$ relative to $A$ has been considered \cite{Boca, Popa06} and proved to be very useful.
It is in particular one of the crucial ingredients used by Popa \cite{Popa06}, to provide the first example of a $II_1$ factor with trivial fundamental group.

A countable measured groupoid $(G,\mu)$ with set of units $X$ (see section \ref{prelim}) gives rise to an inclusion $A\subset M$, 
where $A = L^\infty(X,\mu)$ and where $M = L(G,\mu)$ is the von Neumann algebra of the groupoid. This inclusion is canonically equipped with
a conditional expectation $E_A : M \to A$. Although $M$ is not always a finite von Neumann algebra, there is still a notion of property (H)
relative to $A$ and $E_A$ (see \cite{Ueda}). However, to our knowledge, this property has not been translated in terms  only involving $(G,\mu)$, as in the group case. A significant exception  concerns the case where $G = \cR$
is a countable equivalence relation on $X$, preserving the probability measure $\mu$ \cite{Joli}. Property (H) has also been considered in the different 
topological setting of locally compact groupoids \cite{Tu, G-K}.

Our first goal is to extend the work of Jolissaint \cite{Joli} in order to cover the general case of countable measured groupoids, and in particular the case of group actions
leaving  a probability measure quasi-invariant. Although it is not difficult to guess the right definition of property (H) for $(G,\mu)$ (see definition \ref{groupoid_Haag}),
it is more intricate to prove the equivalence of this notion with the fact that $L(G,\mu)$ has property (H) relative to $L^\infty(X)$.

First, in section \ref{1}, we introduce the basic notions and notation concerning countable measured groupoids. In particular we discuss the Tomita-Takesaki
theory for their von Neumann algebras. This is essentially a reformulation of the pioneering results of P. Hahn \cite{Hahn} in a way
that fits better for our purpose. Sections \ref{2} and \ref{3} discuss in detail several facts about the von Neumann algebra of the Jones' basic construction
for an inclusion $A\subset M$ of von Neumann algebras, in a general setting  and then in our situation $L^\infty(X,\mu) \subset L(G,\mu)$. These are mainly
folklore results but not always easy to find in the literature, and we give a more handy exposition for our purposes. 

In sections \ref{4} and \ref{5}, we study the relations between positive definite functions on our groupoids and completely positive maps on the
corresponding von Neumann algebras. These results are extensions of well known results for groups and of results obtained by
Jolissaint in \cite{Joli} for equivalence relations, but additional difficulties must be overcome. After this preliminary work, it is immediate (section \ref{6})
to show the equivalence of our definition of property (H) for groupoids with the definition involving operator algebras (theorem \ref{equiv_def_property_H}).

Our main motivation originates from the reading of  Ueda's paper \cite{Ueda} and concerns treeable groupoids. This notion was
introduced by Adams for probability measure preserving countable equivalence relations \cite{Ad}. Treeable groupoids may be viewed as the groupoid analogue
of free groups. So a  natural question, raised by C.C. Moore in his survey \cite[p. 277]{Moore08} is whether a treeable equivalence relation must have the Haagerup property.
In fact, this problem is solved in \cite{Ueda} using operator algebras techniques. In Ueda's paper, the notion of treeing is translated in an operator algebra
framework regarding the inclusion $L^\infty(X,\mu) \subset L(G,\mu)$, and it is proved that this condition implies that $L(G,\mu)$ has the Haagerup property
relative to $L^\infty(X,\mu)$. 

Our approach is opposite. For us, it seems more natural to   compare these two notions, treeability and property (H), purely at the level of the groupoid.
Indeed, the definition of treeability is more nicely read at the level of the groupoid than at the level of its von Neumann algebra : roughly speaking, it means that there is a measurable way to endow each fibre of the groupoid with a structure of tree (see definition \ref{def_treeable}). The direct proof that treeability implies property (H)
is given in section \ref{8} (theorem \ref{Ueda}).

For the sake of completeness, we provide the expected results comparing amenabi\-lity, property (H) and property (T) of Kazhdan.
Although the notion amenable groupoid has been extensively studied in \cite{AD-R}, the needed characterization allowing to show that an amenable
countable measured groupoid has the Haagerup property is missing there. We fill this gap in section \ref{7}. We point out that using operator algebras
techniques, this implication is not obvious. Indeed, the amenability of $(G,\mu)$ implies that $L(G,\mu)$ is amenable and so 
amenable relative to $L^\infty(X,\mu)$ (see \cite[\S 3.2]{Popa_Cor}). But in general, this notion is too weak to imply the relative Haagerup property.
What is needed is the notion of strong relative amenability introduced by Popa. We refer the reader to   \cite[rem. 3.5]{Popa06} for more comments
on this fact.

Using our previous work on groupoids with property (T),  we prove in section \ref{9} that, under an assumption of ergodicity, this property is incompatible
with the Haagerup property (theorem \ref{H-T}). As a consequence, we recover the result of Jolissaint \cite[prop. 3.2]{Joli} stating that if $\Gamma$ is a Kazhdan countable group
which acts ergodically on a Lebesgue space $(X,\mu)$ and leaves the probability measure $\mu$ invariant, then the orbit equivalence relation $(\cR_\Gamma,\mu)$
has not the Haagerup property (corollary \ref{T-non-H}). A fortiori, $(\cR_\Gamma,\mu)$ is not treeable, a result due to Adams and Spatzier \cite[thm. 18]{AS}
and recovered in a different way by Ueda.

A groupoid brings an equivalence relation and a bundle of groups into play. More precisely, it is an extension of its associated equivalence relation
by the bundle of its isotropy groups. Section \ref{10} contains some observations relating amenability, treeability,  property (H) and property (T) to
the corresponding properties for the associated equivalence relation and the isotropy groups.

For the reader's convenience, we have gathered in an appendix the needed definitions concerning $G$-bundles of graphs and representations of groupoids,.
 
\setcounter{equation}{0}
\renewcommand{\theequation}{\thesection.\arabic{equation}}

\section{The von Neumann algebra of a measured groupoid}\label{1}

\subsection{Preliminaries on countable measured groupoids} \label{prelim}
Our references for measured groupoids are \cite{AD-R, Hahn, Ram}. Let us first introduce some notation. 
Given a groupoid
$G$, $G^{(0)}$ denotes its unit space and $G^{(2)}$ the set of composable pairs.
The range and source maps from $G$ to  $G^{(0)}$ are denoted respectively
by $r$ and $s$. The corresponding fibres are denoted respectively by $G^x =
r^{-1}(x)$ and $G_x = s^{-1}(x)$. Given subsets $A, B$ of $G^{(0)}$, we define $G^A =
r^{-1}(A)$,
$G_B = s^{-1}(B)$ and $G_{B}^A = G^A \cap G_B$. For $x\in G^{(0)}$, the isotropy
group $G^{x}_x$ is denoted by $G(x)$. The
{\it reduction} of $G$ to
$A$ is the groupoid $G|_A = G_{A}^A$. Two elements $x,y$ of $G^{(0)}$ are said to be
equivalent if
$G^{x}_y \not = \emptyset$. We denote by $\cR_G$  this equivalence relation. If $A \subset G^{(0)}$, its {\it saturation} $[A]$ is the
set $s(r^{-1}(A))$ of
all elements in $G^{(0)}$ that are equivalent to some element of $A$. When $A = [A]$,
we say that $A$ is {\it invariant}.

A Borel groupoid is a groupoid $G$ endowed with a standard Borel structure such that the range, source, inverse and product are Borel maps,
where $G^{(2)}$ has the Borel strutrure induced by $G\times G$ and $G^{(0)}$ has the Borel structure induced by $G$.
We say that $G$ is {\it countable} (or {\it discrete}) if the fibres $G^x$ (or equivalently $G_x$) are countable. 

In the sequel, we only consider such groupoids. We {\it always denote by $X$ the set $G^{(0)}$ of units of $G$}.
A {\it bisection} $S$ is a Borel subset of $G$ such that the restrictions of $r$ and $s$ to $S$ are injective. 
A useful fact, consequence of a theorem of Lusin-Novikov states that, since $r$ and $s$ are countable-to-one Borel maps 
between standard Borel spaces, there exists a countable partition of $G$ into bisections (see \cite[thm. 18.10]{Kech}).

Let $\mu$ be a probability measure on  on $X = G^{(0)}$. We define a $\sigma$-finite measure $\nu$ on $G$ by the formula
$$\int_G F \rd \nu = \int_X \big(\sum_{s(\gamma) = x} F(\gamma)\big) \rd \mu(x).$$

We say that $\mu$ is quasi-invariant if $\nu$ is equivalent to it image $\nu^{-1}$ under $\gamma \mapsto \gamma^{-1}$. 
This is equivalent to the fact that for every bisection $S$, one has $\mu(s(S)) = 0$ if and only if $\mu(r(S)) = 0$.
We set $\delta = \frac{\rd \nu^{-1}}{\rd \nu}$.

\begin{defn} A {\it  countable (or discrete) measured groupoid}\footnote{In \cite{AD-R}, a countable  measured groupoid is called $r$-discrete.
Another difference is that we have swapped here the definitions of $\nu$ and $\nu^{-1}$.} $(G,\mu)$ is a Borel groupoid $G$ equipped with a quasi-invariant probability measure $\mu$ on $X = G^{(0)}$.
\end{defn}

Note that if $\mu_1$ is a probability measure on $X$ with density  $h$ with respect to  $\mu$, then the corresponding measure $\nu_1$ on $G$
has density  $h\circ s$ with respect to $\nu$. In particular, quasi-invariance only depends on the measure class of $\mu$.

\begin{exs} (a) Let $\Gamma \actson X$ be a (right) action of a countable group $\Gamma$ on a standard Borel space $X$, and  assume that the action preserves the class of a probability measure $\mu$. Let $G = X\rtimes \Gamma$ be
the {\it semi-direct product groupoid}. We have $r(x,t) = x$ and $s(x,t) = xt$. The product is given by the formula
$(x,s)(xs, t) = (x, st)$. Equipped with the quasi-invariant measure $\mu$, $(G,\mu)$ is a countable measured groupoid. 
We have $\delta(x,t) =\big( \rd(t^{-1}\mu)/\rd\mu\big)(x)$. As a particular case, we find the group $G = \Gamma$.

(b) Another important family of examples concerns the {\it countable measured equivalence relations}. We are given
an equivalence relation $\cR\subset X\times X$ on a standard Borel space $X$, 
which is a Borel subset of $X\times X$ and whose
equivalence classes are finite or countable. It has an obvious structure of Borel groupoid with $r(x,y) = x$, $s(x,y) = y$
and $(x,y)(x,z) = (x,z)$. When equipped with a quasi-invariant probability measure $\mu$, we say that $(\cR,\mu)$
is a countable measured equivalence relation. Here, quasi-invariance also means that for every Borel subset $A\subset X$, we have $\mu(A) = 0$ if and only if the measure of the saturation $s(r^{-1}(A))$ of $A$ is still $0$.
\end{exs}

A general groupoid is a combination of an equivalence relation with groups.  
Indeed, let $(G,\mu)$ be a countable measured groupoid. Let $c = (r,s)$ be the map $\gamma \mapsto (r(\gamma), s(\gamma))$
from $G$ into $X\times X$. The range of $c$ is the graph $\cR_G$ of the equivalence relation  induced on $X$ by $G$.
Moreover $(\cR_G,\mu)$ is a countable measured equivalence relation. The kernel of the groupoid homomorphism $c$
is the isotropy bundle $(G(x))_{x\in X}$. So, $(G,\mu)$ appears as the extension of the equivalence relation $(\cR,\mu)$ by the isotropy bundle.

 A reduction $(G_{|_U}, \mu_{|_U})$ such that $U$ is conull in $X$ is called {\it inessential}. Since we are working in the setting of measured spaces,
 it will make no difference to replace $(G,\mu)$ by any of its inessential reductions.
 
\subsection{The von Neumann algebra of $(G,\mu)$} If $f:  G\to \C$ is a Borel function, we set
$$\norm{f}_I = \max \set{\norm{x \mapsto \sum_{r(\gamma) = x} \abs{f(\gamma)}}_\infty,\,\, \norm{x \mapsto \sum_{s(\gamma) = x} \abs{f(\gamma)}}_\infty}.$$
Let $I(G)$ be the set of functions such that $\norm{f}_I < +\infty$. It only depends on the measure class of $\mu$.
We endow $I(G)$ with the (associative) convolution product
$$(f*g)(\gamma) = \sum_{\gamma_1\gamma_2 = \gamma} f(\gamma_1)g(\gamma_2) = \sum_{s(\gamma) = s(\gamma_2)} f(\gamma\gamma_{2}^{-1}) g(\gamma_2)
= \sum_{r(\gamma_1) = r(\gamma)} f(\gamma_1) g(\gamma_{1}^{-1} \gamma).$$
and the involution
$$f^*(\gamma) = \overline{f(\gamma^{-1})}.$$
 
For $f$, $g\in I(G)$ we have
\begin{align*}
\sum_{r(\gamma) = x} \abs{f*g)(\gamma)} & \leq \sum_{\set{\gamma,\gamma_1 :\, r(\gamma) = r(\gamma_1) = x}} \abs{f(\gamma_1)}\abs{g(\gamma_1^{-1}\gamma)}\\
&\leq \sum_{r(\gamma_1) = x} \abs{f(\gamma_1)} \sum_{r(\gamma) = x} \abs{g(\gamma_1^{-1}\gamma)} \leq \norm{f}_I\norm{g}_I ,
\end{align*}
and similarly $\sum_{s(\gamma) = x} \abs{f*g)(\gamma)}\leq \norm{f}_I\norm{g}_I$, whence $\norm{f*g}_I \leq \norm{f}_I \norm{g}_I$.

We have $I(G)\subset L^1(G,\nu)\cap L^\infty(G,\nu) \subset L^2(G,\nu)$, with $\norm{f}_1 \leq \norm{f}_I$. Therefore $\norm{\cdot}_I$
is a norm on $I(G)$, where two functions which coincide $\nu$-almost everywhere are identified. Let us show that $I(G)$ is complete for the norm $\norm{\cdot}_I$.
Let $(f_n)$ be a Cauchy sequence. Then $f_n \to f$ and $f_{n}^* \to g$ in $L^1(G,\nu)$. Taking subsequences which converges almost everywhere, we see
that $g = f^*$ almost everywhere. It is a routine exercise to check that $f\in I(G)$ and that $\lim_n \norm{f_n - f}_I = 0$.

Therefore $(I(G), \norm{\cdot}_I)$ is a Banach $*$-algebra. This variant of the Banach algebra $I(G)$ introduced by Hahn \cite{Hahn} has been considered
by Renault in \cite[p. 50]{Ren}. Its advantage is that it does not involve the Radon-Nikodym derivative $\delta$.

For $f\in I(G)$ and $\xi \in L^2(G,\nu)$ we set 
\begin{equation}\label{reg_rep}
(L(f) \xi)(\gamma) = (f*\xi)(\gamma)=  \sum_{\gamma_1\gamma_2 = \gamma} f(\gamma_1)\xi(\gamma_2).
\end{equation}
This defines a bounded operator on $L^2(G,\nu)$. Indeed, given $\eta \in L^2(G,\nu)$
 we have\footnote{When there is no risk of ambiguity, we denote by $\norm{\cdot}_2$ the norm in $L^2(G,\nu)$.}

 \begin{align*}
& \scal{\abs{\eta}, L(\abs{f}) \abs{\xi}}  =\int_X\Big(\sum_{\gamma_1, \gamma_2\in G_x} \abs{\eta(\gamma_1)}\abs{\xi(\gamma_2)} \abs{f(\gamma_1\gamma_2^{-1})}\Big) \rd\mu(x)\\
 &\leq \int_X\Big(\sum_{\gamma_1, \gamma_2\in G_x} \abs{\eta(\gamma_1)}^2\abs{f(\gamma_1\gamma_2^{-1})}\Big)^{1/2}
 \Big(\sum_{\gamma_1, \gamma_2\in G_x} \abs{\xi(\gamma_2)}^2 \abs{f(\gamma_1\gamma_2^{-1})}\Big)^{1/2}\rd\mu(x)\\
  &\leq \Big(\int_X\sum_{\gamma_1, \gamma_2\in G_x} \abs{\eta(\gamma_1)}^2 \abs{f(\gamma_1\gamma_2^{-1})}\rd\mu(x)\Big)^{1/2}
\Big(\int_X \sum_{\gamma_1, \gamma_2\in G_x} \abs{\xi(\gamma_2)}^2 \abs{f(\gamma_1\gamma_2^{-1})}\rd\mu(x)\Big)^{1/2}\\
& \leq \norm{\eta}_2\norm{f}_I^{1/2} \norm{\xi}_{2} \norm{f}_I^{1/2}.
 \end{align*} 
 
 Hence $\norm{L(f)} \leq \norm{f}_I$. 
 One has $L(f)^* = L(f^*)$ ans $L(f) L(g) = L(f*g)$. Hence, $L$ is a representation of $I(G)$, called the {\it left regular representation}.
 
 Note that $L^2(G,\nu)$ is a direct integral of Hilbert spaces :
 $$L^2(G,\nu) = \int_X^\oplus \ell^2(G_x) \rd\mu(x).$$
 Under the product
 $$(f\xi)(\gamma) = f\circ s(\gamma) \xi(\gamma)$$
 where $f\in L^\infty(X)$ and $\xi\in L^2(G,\nu)$, we define on $L^2(G,\nu)$ a structure of $L^\infty(X)$-module.
 Obviously, the representation $L$ commutes with this action of $L^\infty(X)$. In fact $L^\infty(X)$ is the algebra of diagonalizable
 operators with respect to the disintegration $L^2(G,\nu) = \int_X^\oplus \ell^2(G_x) \rd\mu(x)$.
 
 \begin{defn} The {\it von Neumann algebra of the countable measured groupoid} $(G,\mu)$ is the von Neumann subalgebra $L(G,\mu)$ of $\cB(L^2(G,\nu))$
 generated by $L(I(G))$. It will also be denoted by $M$ in the rest of the paper.
 \end{defn}
 
 Since the elements of $L(G,\mu)$ commute with the diagonal action of $L^\infty(X)$, they are decomposable operators (\cite[thm 1, p. 164]{Dix}).
 In fact, $L(f) = \int_X^{\oplus} L_x(f) \rd\mu(x)$, where $L_x(f): \ell^2(G_x) \to \ell^2(G_x)$ is defined as in \eqref{reg_rep}, but for $\xi\in \ell^2(G_x)$.

Let $C_n = \set{1/n \leq \delta \leq n}$. Then $(C_n)$ is an increasing sequence of measurable subsets of $G$ with $\cup_n C_n = G$ (up
to null sets). We denote by $I_n(G)$ the set of elements in $I(G)$ taking value $0$ outside $C_n$ and we set $I_\infty(G) = \cup_n I_n(G)$.
Obviously, $I_\infty(G)$ is an involutive subalgebra of $I(G)$. We leave it to the reader to check that $I_\infty(G)$
is dense into $L^2(G,\nu)$ and that $L(G,\mu)$ is generated by $L(I_\infty(G))$.

 \begin{rem} \label{Hahn} As already said, our notations differ from that of the founding paper \cite{Hahn} of Hahn. First, we have swapped the roles of $\nu$ and $\nu^{-1}$. Moreover, in \cite{Hahn}, another version of $I(G)$ is used, which however contains $I_\infty(G)$. The regular representation consi\-dered by Hahn, that we denote in this remark by $L_H$ to avoid ambiguity, acts on $L^2(G,\nu^{-1})$
 (with our notation for $\nu$). More precisely, for $f\in I_\infty(G)$, we have
 $$\forall \xi \in L^2(G,\nu^{-1}), \quad L_H(f)\xi = f*\xi \in L^2(G,\nu^{-1}).$$
 One easily transfers the objects and results of \cite{Hahn}
 to ours by using the isometry $V : \xi \mapsto \delta^{1/2} \xi$ from $L^2(G,\nu^{-1})$ to $L^2(G,\nu)$. The interested reader will immediately check
  that for $f\in I_\infty(G)$,
  $$V^* L(f) V = L_H(\delta^{-1/2} f).$$
 Our motivation for our changes of presentation is that, on one hand, we avoid the use of $\delta$ in the definition of $L(G,\mu)$, and that, on the other hand,  the von Neumann algebra $L(G,\mu)$  is made of decomposable operators on $\int_X^\oplus \ell^2(G_x) \rd\mu(x)$.
\end{rem}
 
 The von Neumann algebra $L^\infty(X)$ is isomorphic to a subalgebra of $I_\infty(G)$, by giving to $f\in L^\infty(X)$ the value $0$ outside $X \subset G$.
 Note that, for $\xi\in L^2(G,\nu)$, 
  $$(L(f)\xi)(\gamma) = f\circ r(\gamma) \xi(\gamma).$$
 In this way, $A= L^\infty(X)$ appears as a von Neumann subalgebra of $M$.

 If we replace $\mu$ by an equivalent  probability measure $\mu_1 = h\mu$, we see that $U : \xi \mapsto (h\circ s)\, \xi$ is a unitary operator which carries the convolution product by $f$ in $L^2(G,\nu)$ to the same convolution product in $L^2(G,\nu_1)$. Hence, the
 pair $A\subset M$ only depends on the measure class of $\mu$, up to unitary equivalence.

We  view $I(G)$ as a subspace of $L^2(G,\nu)$. The characteristic function $\un_X$ of $X \subset G$ is a norm one vector in $L^2(G,\nu)$.
 Let $\varphi$ be the normal state on $M$ defined by
 $$\varphi(T) = \scal{\un_X, T \un_X}_{L^2(G,\nu)}.$$
 For $f\in I(G)$, we have
 $$\varphi(L(f)) = \int_X f(x) \rd \mu(x),$$
 and therefore, for $f,g\in I(G)$,
 \begin{equation}\label{scalar}
 \varphi(L(f)^*L(g)) = \scal{f,g}_{L^2(G,\nu)}.
 \end{equation}

\begin{lem}\label{right_action} Let $g$ be a Borel function on $G$ such that $\delta^{-1/2} g = f \in I(G)$ (for instance $g\in I_\infty(G)$). Then $\xi \mapsto \xi * g$ is a bounded operator on $L^2(G,\nu)$. More precisely, we have $$\norm{\xi* g}_2 \leq \norm{f}_I \norm{\xi}_2.$$
\end{lem}
\begin{proof} We have
\begin{align*}
\norm{\xi* g}_{2}^2 & \leq \int_X \sum_{s(\gamma) = x} \Big( \sum_{s(\gamma_1) = x}\abs{\xi(\gamma\gamma_1^{-1})}\abs{f(\gamma_1)} \delta(\gamma_1)^{1/2}\Big)^2 \rd \mu(x)\\
&\leq  \int_X \sum_{s(\gamma) = x} \big( \sum_{s(\gamma_1) = x}\abs{\xi(\gamma\gamma_1^{-1})}^2\abs{f(\gamma_1)} \delta(\gamma_1)\big)
\big(\sum_{s(\gamma_1) = x}\abs{f(\gamma_1)}\big) \rd\mu(x)\\
&\leq \norm{f}_I \int_G\Big(\sum_{\set{\gamma :\, s(\gamma) =s(\gamma_1)}} \abs{\xi(\gamma\gamma_1^{-1})}^2 \Big) \abs{f(\gamma_1)} \rd \nu^{-1}(\gamma_1)\\
&\leq \norm{f}_I \int_X \sum_{s(\gamma_1) = x} \abs{f(\gamma_1^{-1}} \big(\sum_{\set{ \gamma : s(\gamma) = r(\gamma_1)}} \abs{\xi(\gamma\gamma_1)}^2\big) \rd\mu(x)\\
&\leq \norm{f}_I \int_X \sum_{s(\gamma_1) = x} \abs{f(\gamma_1^{-1}} \big(\sum_{s(\gamma_2) = x} \abs{\xi(\gamma_2)}^2\big) \rd\mu(x)\\
&\leq \norm{f}_I^{2} \norm{\xi}_{2}^2.
\end{align*}
\end{proof}

We set $R(g)(\xi) = \xi*g$. We have $L(f)\circ R(g) = R(g) \circ L(f)$
for every $g\in I_\infty(G)$ and $f\in I(G)$.
We denote by $R(G,\mu)$ the von Neumann algebra generated by $R(I_\infty(G))$.

\begin{lem} The vector $\un_X$ is cyclic and separating for $L(G,\mu)$, and therefore $\varphi$ is a faithful state.
\end{lem}

\begin{proof} This follows immediately from the fact that $L(f)$ and $R(g)$ commute for $f,g\in I_\infty(G)$, with
$L(f) \un_X = f$ and $R(g)\un_X = g$, and from the density of $I_\infty(G)$ into $L^2(G,\nu)$.
\end{proof}

The von Neumann algebra $L(G,\mu)$ is on standard form on $L^2(G,\nu)$, canonically identified with $L^2(M,\varphi)$
(see \eqref{scalar}).
We identify $M$ with a dense subspace of $L^2(G,\nu)$   by $T\mapsto \widehat{T} = T(\un_X)$.
The modular conjugation $J$ and the one-parameter modular
group $\sigma$ relative to the vector $\un_X$ (and $\varphi$) have been computed in \cite{Hahn}. With our notations, we have
$$\forall \xi\in L^2(G,\nu), \quad (J\xi)(\gamma) = \delta(\gamma)^{1/2}\overline{\xi(\gamma^{-1})}$$
and 
$$\forall T\in L(G,\mu),\quad \sigma_t(T) = \delta^{it}T\delta^{-it}.$$
 Here, for $t\in \R$, the function $\delta^{it}$ acts on $L^2(G,\nu)$ by pointwise multiplication and defines a unitary operator.
 Note that for $f\in L(G,\mu)$, we  have $\delta^{it} L(f) \delta^{-it} = L(\delta^{it}f)$. 
 In particular,  $\sigma$ acts trivially on $A$. Therefore (see \cite{Tak72}),  there exists a unique faithful conditional expectation $E_A : M\to A$
 such that $\varphi = \varphi\circ E_A$, and 
  for $T\in M$, we have
 $$\widehat{E_A(T)} = e_A(\widehat{T}),$$
where $e_A$ is the orthogonal projection from $L^2(G,\nu)$ onto $L^2(X,\mu)$. If we view the elements of $M$
as functions on $G$, then $E_A$ is the restriction map to $X$. The triple $(M, A, E_A)$ only depends on the class of $\mu$,
up to equivalence.

 For $f\in I(G)$ and $\xi\in L^2(G,\nu)$ we observe that
 $$(JL(f)J)\xi(\gamma) = \sum_{s(\gamma_1) = s(\gamma)} \overline{f(\gamma_1^{-1})} \delta(\gamma_1)^{1/2} \xi(\gamma\gamma_1^{-1}),$$
that is 
\begin{equation}\label{convol_right}
(JL(f)J)\xi = R(g)\xi = \xi * g\quad\hbox{with} \quad g=\delta^{1/2} f^* .
\end{equation}
 
   \section{Basic facts on the module $L^2(M)_A$}\label{2}
   
   We consider, in an abstract setting, the situation we have met above. Let $A\subset M$ be a pair of von Neumann algebras, where $A = L^\infty(X,\mu)$ is abelian.
   We assume the existence of a normal faithful conditional expectation $E_A : M\to A$ and we set $\varphi = \tau_\mu \circ E_A$,
   where $\tau_\mu$ is the state on $A$ defined by the probability measure $\mu$. Recall that $M$ is on standard form on the Hilbert space $L^2(M,\varphi)$
   of the Gelfand-Naimark-Segal construction associated with $\varphi$.
   We view $L^2(M,\varphi)$ as a left $M$-module and a right $A$-module. Identifying\footnote{When necessary, we shall write $\wh{m}$ the element $m\in M$,
  when  viewed in $L^2(M,\varphi)$, in order to stress this fact.} $M$ with a subspace of $L^2(M,\varphi)$, we know that $E_A$ is the restriction to $M$
   of the orthogonal projection $e_A : L^2(M,\varphi) \to L^2(A,\tau_\mu)$.
   
  For further use, we make the following observation
  \begin{equation}\label{droite}
  \forall m\in M, \forall a\in A, \quad JaJ \wh m = \wh m a^* = \wh{ma^*}.
  \end{equation}
  Indeed, if $S$ is the closure of the map $\wh m \mapsto \wh{m^*}$ and if $S = J \Delta^{1/2}= \Delta^{-1/2}J$ is its polar decomposition, then
   every $a\in A$ commutes with $\Delta$ because $\varphi = \tau_\mu \circ E_A$ and $\tau_\mu$ is a trace (see \cite{Tak72}).
  Then \eqref{droite} follows easily. Note that \eqref{convol_right} gives a particular case of this remak.
   
   \subsection{The commutant $\scal{M,e_A}$ of the right action}
   The algebra of all operators which commute with the right action of $A$ is the von Neumann algebra of the basic construction for
  $A\subset M$. It is denoted $\scal{M,e_A}$ since it is generated by $M$ and $e_A$. The linear span of $\set{m_1 e_A m_2 : m_1,m_2\in M}$
  is a  $*$-subalgebra which is weak operator dense in $\scal{M,e_A}$. Moreover $\scal{M,e_A}$ is a semi-finite von Neumann algebra,
  carrying a canonical normal faithful semi-finite trace $\Tr_\mu$ (depending on the choice of $\mu$), defined by
  $$\Tr_\mu(m_1 e_A m_2) = \int_X E_A(m_2 m_1) \rd\mu = \varphi(m_2 m_1).$$
  (for these classical results, see \cite{Jones}, \cite{Popa_book}).  We shall give more information on this trace in lemma \ref{calcul_trace} and its proof.
  We need some preliminaries.
    
    \begin{defn} A vector $\xi\in L^2(M,\varphi)$ is  {\it $A$-bounded} if there exists $c>0$ such that for every $a\in A$,
 $$\norm{\xi a}_2 \leq c \tau_\mu(a^*a)^{1/2}.$$
\end{defn}

(Again, if the norm is clear from the context, we use the notation $\norm{\cdot}_2$ instead of $\norm{\cdot}_{L^2(M)}$.)

We denote by $L^2(M,\varphi)^0$, or $\cL^2(M,\varphi)$, the subspace of $A$-bounded vectors. It contains $M$. We also recall the obvious fact that $T \mapsto T(1_A)$
is an isomorphism from the space $\cB(L^2(A,\tau_\mu)_A, L^2(M,\varphi)_A)$ of $A$-linear bounded operators $ T : L^2(A,\tau_\mu)\to L^2(M,\varphi)$ onto $\cL^2(M,\varphi)$.
For $\xi\in \cL^2(M,\varphi)$, we denote by $L_\xi$ the corresponding operator from $L^2(A,\tau_\mu)$ into $L^2(M,\varphi)$.
  In particular, for $m\in M$, we have $L_{m} = m_{|_{L^2(A,\tau_\mu)}}$. It is easy to see that $\cL^2(M,\varphi)$ is stable under the actions of 
  $\scal{M,e_A}$ and $A$, and that $L_{T\xi a} = T \circ L_\xi \circ a$ for $T\in \scal{M,e_A}$, $\xi \in \cL^2(M,\varphi)$, $a\in A$.

 For $\xi,\eta \in \cL^2(M,\varphi)$, the operator $L_{\xi}^* L_\eta \in \cB(L^2(A,\tau_\mu))$ commutes with $A$ and so is in $A$. We set
 $\scal{\xi, \eta}_A = L_{\xi}^*L_\eta$. In particular, we have $\scal{{m_1}, {m_2}}_A = E_A(m_{1}^*m_2)$ for $m_1,m_2\in M$.
   Equipped with the $A$-valued inner product
 $$\scal{\xi, \eta}_A = L_{\xi}^*L_\eta,$$
  $\cL^2(M,\varphi)$ is a self-dual Hilbert right $A$-module. It is a normed space with respect to the norm
  $$\norm{\xi}_{\cL^2(M)} = \norm{\scal{\xi,\xi}_A}_{A}^{1/2}.$$ 
  Note that 
  \begin{equation}\label{compr_norm}
  \norm{\xi}_{L^2(M)} = \norm{\scal{\xi,\xi}_A}_{L^2(A)}^{1/2} \leq   \norm{\xi}_{\cL^2(M)}.
  \end{equation}
  
  We define a positive hermitian form on the algebraic tensor product $\cL^2(M,\varphi)\odot L^2(A)$ by
  $$\scal{\xi\otimes f, \eta\otimes g} =  \int_X \overline{f} g \scal{\xi,\eta}_A \rd \mu .$$
  The Hilbert space $\cL^2(M,\varphi)\otimes_A L^2(A)$ obtained by separation and completion is isomorphic to $L^2(M,\varphi)$ as a right $A$-module by
  $$\xi\otimes {f} \mapsto \xi f.$$ Moreover the von Neumann algebra $\cB(\cL^2(M,\varphi)_A)$ of bounded $A$-linear
  endomorphisms of $\cL^2(M,\varphi)$ is isomorphic to $\scal{M,e_A}$ by
  $$T \mapsto T\otimes 1.$$
We shall identify these two von Neumann algebras (see \cite{Pas73}, \cite{Rie74} for details on these facts).

\begin{defn} An {\it orthonormal basis} of the $A$-module $L^2(M,\varphi)$ is a family $(\xi_i)$ of elements of $\cL^2(M,\varphi)$ 
such that $\overline{\sum_i\xi_i A} = L^2(M,\varphi)$ and $\scal{\xi_i,\xi_j}_A =\delta_{i,j} p_j$ for all $i,j$, where the $p_j$ are projections in $A$.
\end{defn}

It is easily checked that $L_{\xi_i}L_{\xi_i}^*$ is the orthogonal projection on $\overline{\xi_i A}$, and that these projections are mutually orthogonal
with $\sum_i L_{\xi_i}L_{\xi_i}^* = 1$.

Using a generalization of the Gram-Schmidt orthonormalization process, one shows the existence of orthonormal bases (see \cite{Pas73}).

\begin{lem}\label{calcul_scal} Let $(\xi_i)$ be an orthonormal basis of the $A$-module $L^2(M,\varphi)$. For every $\xi\in \cL^2(M,\varphi)$, we have
\begin{equation}\label{eq_calcul_scal}
\scal{\xi,\xi}_A = \sum_i \scal{\xi,\xi_i}_A \scal{\xi_i, \xi}_A
\end{equation}
(weak* convergence).
\end{lem}

\begin{proof} We have $\sum_i L_{\xi_i} L_{\xi_i}^* = 1$ and 
$$\scal{\xi,\xi}_A = L_{\xi}^* L_{\xi} = \sum_i (L_{\xi}^* L_{\xi_i}) (L_{\xi_i}^*L_\xi) = \sum_i \scal{\xi,\xi_i}_A \scal{\xi_i, \xi}_A.$$
\end{proof}

\begin{lem}\label{calcul_trace} Let $(\xi_i)$ be an orthonormal basis of the $A$-module $L^2(M,\varphi)$. 
\begin{itemize}
\item[(i)] For every $x\in \scal{M,e_A}_+$
we have
\begin{equation}\label{eq_calcul_trace1}
\Tr_\mu(x) = \sum_i \tau_\mu(\scal{\xi_i, x\xi_i}_A) = \sum_i \scal{\xi_i, x\xi_i}_{L^2(M)}.
\end{equation}
\item[(ii)] $\sp\set{L_\xi L_{\eta}^*: \xi,\eta \in \cL^2(M,\varphi)}$ is contained in the ideal of definition of $\Tr_\mu$
and we have, for $\xi,\eta \in  \cL^2(M,\varphi)$,
\begin{equation}\label{eq_calcul_trace2}
\Tr_\mu(L_\xi L_{\eta}^*) = \tau_\mu( L_\eta^{*} L_\xi) = \tau_\mu(\scal{\eta,\xi}_A).
\end{equation}
\end{itemize}
\end{lem}

\begin{proof} (i) The map $U : L^2(M,\varphi) = \oplus_i \overline{\xi_i A} \to \oplus_i p_i L^2(A)$ defined by $U(\xi_i a) = p_i a$
is an isomorphism which identifies $L^2(M,\varphi)$ to the submodule $p(\ell^2(I) \otimes L^2(A))$ of $\ell^2(I) \otimes L^2(A)$, with $p = \oplus_i p_i$.
The canonical trace on $\scal{M,e_A}$ is transfered to the restriction to $p\big(\cB(\ell^2(I)\ovtimes A\big)p$ of the trace $\Tr\otimes \tau_\mu$,
defined on $T = [T_{i,j}] \in (\cB(\ell^2(I)\ovtimes A)_+$ by
 $$(\Tr\otimes \tau_\mu)(T) = \sum_i \tau_\mu(T_{ii}).$$
 It follows that
 $$\Tr_\mu(x) = \sum_i \tau_\mu((UxU^*)_{ii}) =  \sum_i \scal{\xi_i, x\xi_i}_{L^2(M)} =  \sum_i \tau_\mu(\scal{\xi_i, x\xi_i}_A).$$
 
 (ii) Taking $x = L_\xi L_\xi^{*}$ in (i), the equality $\Tr_\mu(L_\xi L_{\xi}^*) = \tau_\mu(\scal{\xi,\xi}_A)$
 follows from equations \eqref{eq_calcul_scal} and \eqref{eq_calcul_trace1}. Formula \eqref{eq_calcul_trace2}
 is deduced by polarization.
 \end{proof}
 
\subsection{Compact operators}
In a semi-finite von Neumann algebra $N$, there is a natural notion of ideal of compact operators,
namely the  norm-closed ideal $\cI(N)$ generated by its  finite projections (see \cite[\S 1.3.2]{Popa06} or
\cite{P-R}).  

For $N = \scal{M, e_A}$, there is another natural candidate for the space of compact operators.
First, we observe that given $\xi,\eta\in \cL^2(M,\varphi)$, the operator $L_\xi L_{\eta}^*\in \scal{M, e_A}$ plays the role of a rank one operator in ordinary
Hilbert spaces : indeed, if $\alpha \in \cL^2(M,\varphi)$, we have $(L_\xi L_{\eta}^*)(\alpha) = \xi \scal{\eta,\alpha}_A$.
In particular, for $m_1,m_2\in M$, we note that $m_1 e_A m_2$ is a ``rank one operator'' since $m_1 e_A m_2 = L_{m_1}L_{m_{2}^*}^*$.
We denote by $\cK(\scal{M,e_A})$ the norm closure into $\scal{M,e_A}$ of 
$$\sp\set{L_\xi L_{\eta}^*: \xi,\eta \in \cL^2(M,\varphi)}.$$
 It is a two-sided ideal of $\scal{M,e_A}$.

For every $\xi \in \cL^2(M,\varphi)$, we have $L_\xi e_A \in \scal{M,e_A}$. Since
$$L_\xi L_{\eta}^* = (L_\xi e_A) (L_\eta e_A)^*$$
we see that $\cK(\scal{M,e_A})$ is the norm closed two-sided ideal generated by $e_A$ in $\scal{M,e_A}$.
The projection $e_A$ being finite  (because $\Tr_\mu(e_A) = 1)$, we have 
$$\cK(\scal{M,e_A}) \subset \cI(\scal{M,e_A}).$$

The subtle difference between $\cK(\scal{M,e_A})$ and $\cI(\scal{M,e_A})$ is studied in \cite[\S 1.3.2]{Popa06}.  We recall
in particular that for every $T \in \cI(\scal{M,e_A})$ and every $\varepsilon >0$, there is a projection $p\in A$
such that $\tau_\mu(1- p)\leq \varepsilon$ and $ TJpJ \in \cK(\scal{M,e_A})$ (see \cite[Prop. 1.3.3 (3)]{Popa06}).
\footnote{In \cite{Popa06}, $\cK(\scal{M,e_A})$ is denoted $\cI_0(\scal{M,e_A}$.}

\subsection{The relative Haagerup property}
Let $\Phi$ be a unital completely positive map from $M$ into $M$ such that $E_A\circ \Phi = E_A$. Then for $m\in M$, we have
$$\norm{\Phi(m)}_2^{2} = \varphi(\Phi(m)^*\Phi(m)) \leq \varphi(\Phi(m^*m)) = \varphi(m^*m) = \norm{m}_2^{2}.$$
It follows that $\Phi$ extends to a contraction $\wh\Phi$ of $L^2(M,\varphi)$. Assume now that $\Phi$ is $A$-bilinear.
Then $\wh\Phi$ commutes with the right action of $A$ (due to \eqref{droite}) and so belongs to $\scal{M,e_A}$. It
also commutes with the left action of $A$ and so belongs to  $A' \cap \scal{M,e_A}$.

\begin{defn}\label{def_hag_rel} We say that $M$ has the {\it Haagerup property} ({\it or property} (H)) {\it relative to $A$ and $E_A$} if there exists a net
$(\Phi_i)$ of unital $A$-bilinear completely positive maps from $M$ to $M$ such that
\begin{itemize}
\item[(i)] $E_A \circ \Phi_i = E_A$ for all $i$ ;
\item[(ii)]  $\widehat{\Phi_i} \in \cK(\scal{M,e_A})$ for all $i$ ;
\item[(iii)] $\lim_i \norm{\Phi_i(x) - x}_2 = 0$ for every $x\in M$.
\end{itemize}
\end{defn}

This notion is due to Boca \cite{Boca}. In \cite{Popa06}, Popa uses a slightly different formulation. 

\begin{lem}\label{equiv_Haag} In the previous definition, we may equivalently assume that $\widehat{\Phi_i} \in \cI(\scal{M,e_A})$ for every $i$.
\end{lem}

\begin{proof} This fact is explained in \cite{Popa06}. Let $\Phi$ be a unital  $A$-bilinear completely positive map from $M$ to $M$  such that $E_A\circ \Phi = E_A$ 
and $\widehat{\Phi} \in \cI(\scal{M,e_A})$. As already said, by \cite[Prop. 1.3.3 (3)]{Popa06}, for every $\varepsilon >0$, there is a projection $p$ in $A$
with $\tau_\mu(1 - p) <\varepsilon$ and $\wh{\Phi}JpJ \in \cK(\scal{M,e_A})$. Thus we have $p\wh{\Phi}JpJ \in \cK(\scal{M,e_A})$. Moreover, this operator is associated with the completely positive map $\Phi_p : m\in M \mapsto \Phi(pmp)$, since
$$(p\wh{\Phi} JpJ)(\wh m) = p\wh{\Phi}(\wh m p) = p\wh{\Phi}(\wh{mp}) = p\wh{\Phi(mp)} = \wh{\Phi(pmp)}.$$
Then, $\Phi' = \Phi_p + (1-p) E_A$ is unital, satisfies $E_A\circ \Phi' = E_A$
and still provides an element of $\cK(\scal{M,e_A})$. This modification allows to prove that if definition \ref{def_hag_rel} holds
with $\cK(\scal{M,e_A})$ replaced by $\cI(\scal{M,e_A})$, then the relative Haagerup property is satisfied (see \cite[Prop. 2.2 (1)]{Popa06}).
\end{proof} 

     \section{Back to $L^2(G,\nu)_A$}\label{3}
  We apply the  facts just reminded to $M = L(G,\mu)$, which is on standard form on $L^2(G,\nu) = L^2(M,\varphi)$. This Hilbert space
  is viewed as a right $A$-module : for $\xi\in L^2(G,\mu)$ and $f\in A$, the action is given by $\xi f\circ s$.
  
  It is easily seen that $\cL^2(M,\varphi)$ is the space of $\xi\in L^2(G,\nu)$ such that 
  $$x \mapsto \sum_{s(\gamma) = x} \abs{\xi(\gamma)}^2$$
  is in $L^\infty(X)$. Moreover, for $\xi,\eta\in \cL^2(M,\varphi)$ we have
  $$\scal{\xi, \eta}_A = \sum_{s(\gamma) = x} \overline{\xi(\gamma)}\eta(\gamma).$$
  
  For simplicity of notation, we shall often identify $f\in I(G)\subset L^2(G,\nu)$ with the operator $L(f)$.\footnote{The reader should not confuse $L(f) : L^2(G,\nu) \to L^2(G,\nu)$
  with its restriction $L_f : L^2(A,\tau_\mu) \to L^2(G,\nu)$.} For instance, for $f,g\in I(G)$, the operator $L(f)\circ L(g)$
  is also written $f*g$, and for $T\in \cB(L^2(G,\mu))$, we write $T\circ f$ instead of $T\circ L(f)$.
 
 Let $S\subset G$ be a bisection. Its characteristic function $\un_S$ is an element of $I(G)$ and a partial isometry in $M$ since
  \begin{equation}\label{bis_eq}\un_{S}^*  * \un_{S} = \un_{s(S)} ,\quad\hbox{and} \quad \un_{S} * \un_{S}^*  = \un_{r(S)}.
  \end{equation}
  Another straightforward computation shows that $\un_S \circ e_A\circ  \un_{S}^*$ is the projection given by the pointwise multiplication by $\un_S$:
$$\forall \xi \in L^2(G,\nu), \,( (\un_S\circ e_A \circ \un_{S}^*) \xi)(\gamma) = \un_S(\gamma)\xi(\gamma).$$
 This multiplication operator by $\un_S$ is denoted $\mathfrak{m}(\un_S)$.

 Let $G = \sqcup S_n$ be  a countable partition of $G$ into Borel bisections. 
Since $L^2(M,\varphi)$ is the orthogonal direct sum of the subspaces $\overline{\un_{S_n} A}$ and since 
the ``rank one operator'' $\un_{S_n}\circ e_A \circ \un_{S_n}^*$ is the orthogonal projection on $\overline{\un_{S_n} A}$ for
every $n$, we see that $(\un_{S_n})_n$ is an orthonormal basis of the right $A$-module $L^2(M,\varphi)$.

By lemma \ref{calcul_trace}, for $x\in \scal{M,e_A}_+$ we have
$$\Tr_\mu(x) = \sum_n \scal{\un_{S_n}, x\un_{S_n}}_{L^2(M)}.$$
In particular, whenever $x$ is the multiplication operator $\mathfrak{m}(f)$ by some non-negative Borel function $f$, we get
\begin{equation}\label{formule_trace}
\Tr_\mu(\mathfrak{m}(f)) = \int_G f \rd \nu.
\end{equation}

\section{From completely positive maps to positive definite functions}\label{4}
Recall that if $G$ is a countable group, and $\Phi : L(G) \to L(G)$ is a completely positive map, then $t\mapsto F_\Phi(t)= \tau(\Phi(u_t) u_{t}^*)$
is a positive definite function on $G$, where $\tau$ is the canonical trace on $L(G)$ and $u_t$, $t\in G$, are the canonical unitaries.
We want to extend this classical fact to the groupoid case. This was achieved by Jolissaint \cite{Joli} for countable probability measure preserving equivalence relations.

Let $(G,\mu)$ be a countable measured groupoid and $M = L(G,\mu)$.
Let $\Phi : M\to M$ be a normal $A$-bilinear unital completely positive map. Let $G = \sqcup S_n$ be a partition into Borel bisections. We define $F_\Phi : G \to \C$ by
\begin{equation}\label{def_fphi}
F_\Phi(\gamma) = E_A(\Phi(\un_{S_n})\circ \un_{S_n}^*)\circ r(\gamma),
\end{equation}
where $S_n$ is the bisection which contains $\gamma$.

That $F_\Phi$ does not depend (up to null sets) on the choice of the partition is a consequence of the following lemma.

\begin{lem}\label{prelim_1} Let $S_1$ and $S_2$ be two Borel bisections. Then
 $$E_A(\Phi(\un_{S_1})\circ\un_{S_1}^*) = E_A(\Phi(\un_{S_2})\circ\un_{S_2}^*)$$
 almost everywhere on $r(S_1\cap S_2)$.
 \end{lem}
 
 \begin{proof} Denote by $e$ the characteristic function of $r(S_1\cap S_2)$. Then
 $e * \un_{S_1} = e*\un_{S_2} = \un_{S_1\cap S_2}$. Thus we have
 \begin{align*}
 eE_A(\Phi(\un_{S_1})\circ\un_{S_1}^*)e & = E_A\big(\Phi(e* \un_{S_1})\circ (\un_{S_1}^* * e)\big)\\
 &= E_A\big(\Phi(e* \un_{S_2})\circ (\un_{S_2}^* * e) =  eE_A(\Phi(\un_{S_2})\circ\un_{S_2}^*)e .
 \end{align*}
  \end{proof}
  
 We now want to show that $F_\Phi$ is a positive definite  function in the following sense. We shall need some preliminary facts.

 \begin{defn} A Borel function $F : G \to \C$ is said to be   {\it positive definite} if  there exists a $\mu$-null subset $N$ of $X= G^{(0)}$ such that for every $x\notin N$,
 and every $\gamma_1,\dots,\gamma_k \in G^x$, the $k\times k$ matrix $[F(\gamma_{i}^{-1}\gamma_j)]$ is positive.
 \end{defn}
 
 \begin{defn}\label{admis} We say that a {\it Borel bisection $S$ is admissible} if there exists a constant $c>0$ such that
 $1/c \leq \delta(\gamma) \leq c$ almost everywhere on $S$.
 \end{defn}
 In other terms, $\un_S \in I_\infty(G)$ and so the convolution to the right by $\un_S$ defines a bounded operator $R(\un_S)$ on $L^2(M,\varphi)$, by \eqref{convol_right}.

 \begin{lem}\label{bon_bis_1} Let $S$ be a Borel bisection and let $T\in M$. We have $\widehat{\un_S \circ T} = \un_S * \widehat{T}$. Moreover, if $S$  is admissible, we have
 $\widehat{T\circ \un_S} = \widehat{T}* \un_S$.
 \end{lem}
 
 \begin{proof} We have $\widehat{\un_S \circ T} = \un_S\circ T(\un_X) = \un_S * \widehat{T}$.
 
 On the other hand, given $f\in I(G)$, we have $\wh{L(f)\circ \un_S} = f*1_S$. So, if $(f_n)$ is a sequence in $I(G)$ such that $\lim_n L(f_n) = T$ in the strong operator topology, we have
 $$\widehat{T\circ \un_S} = T(\wh{\un_S}) = \lim_n L(f_n)(\wh{\un_S})  = \lim_n f_n * \un_S$$
 in $L^2(G,\nu)$. But, when $S$ is admissible, the convolution to the right by $\un_S$ is the bounded operator $R(\un_S)$. Noticing that $\lim_n\norm{f_n - \wh T}_2 = 0$, it follows that 
$$\widehat{T\circ \un_S} = \lim_n f_n * \un_S = \wh T * \un_S.$$ 
  \end{proof}

  \begin{lem}\label{bon_bis_2} Let $T\in M$, and let $S$ be an admissible bisection. Then 
 $$\un_S(\gamma) E_A(T\circ \un_S)(s(\gamma)) = \un_S(\gamma) E_A(\un_S\circ T)(r(\gamma))$$
 for almost every $\gamma$.
 \end{lem}
 
 \begin{proof} We have 
 $$(\wh{T\circ \un_S})(x) = \sum_{\gamma_1\gamma_2 = x} \widehat{T}(\gamma_1) \un_S(\gamma_2)= \widehat{T}(\gamma_2^{-1})$$
 whenever $x\in s(S)$, where $\gamma_2$ is the unique element of $S$ with $s(\gamma_2) = x$. Otherwise 
 $(T\circ \un_S)(x)  = 0$.
 
 On the other hand, 
 $$(\wh{\un_S \circ T})(x) = \widehat{T}(\gamma_1^{-1})$$ whenever $x \in r(S)$, 
 where $\gamma_1$ is the unique element of $S$ with $r(\gamma_1) = x$. Otherwise 
 $(\wh{\un_S \circ T})(x) = 0$. Our statement follows immediately.
 \end{proof}

 \begin{lem} $F_\Phi$ is a positive definite  function.
 \end{lem}
 
 \begin{proof}  We assume that $F_\Phi$ is defined  by equation \eqref{def_fphi} through a partition under admissible bisections.
 We set $S_{ij} = S_i^{-1}S_j = \set{\gamma^{-1}\gamma' : \gamma\in S_{i}, \gamma' \in S_j}$.
 Note that $\un_{S_{i}}^* * \un_{S_j} = \un_{S_{ij}}$. Morever, the $S_{ij}$ are admissible bisections.
 We set 
 $$Z_{ijm} = \set{x\in r(S_{ij}\cap S_{m}) : E_A(\Phi(\un_{S_{ij}})\circ\un_{S_{ij}}^*)(x) \not= E_A(\Phi(\un_{S_{m}})\circ\un_{S_{m}}^*)(x)}$$
 and $Z = \cup_{i,j,m} Z_{ijm}$. It is a null set by lemma \ref{prelim_1}.
 
 By lemma \ref{bon_bis_2}, for every $i$ there is a null set  $E_i \subset r(S_i)$ such that for $\gamma \in S_i$ with $r(\gamma) \notin E_i$ and for every $j$, we have
 $$E_A( \Phi(\un_{S_{ij}}) \circ \un_{S_j}^*\circ \un_{S_i})(s(\gamma)) \not= E_A( \un_{S_i}\circ \Phi(\un_{S_{ij}}) \circ \un_{S_j}^*)(r(\gamma))$$
 We set $E = \cup_i E_i$. Let $Y$ be the saturation of $Z\cup E$. It is a null set, since $\mu$ is quasi-invariant.
  
  Let $x\notin Y$, and $\gamma_1, \dots,\gamma_k\in G^x$.  Assume that $\gamma_{i}^{-1}\gamma_j \in S_{n_i}^{-1}S_{n_j} \cap S_m$.
 We have $r(\gamma_{i}^{-1}\gamma_j) = s(\gamma_i) \notin Y$ since $r(\gamma_i) = x \notin Y$.
 Therefore,
 $$F_{\Phi}(\gamma_{i}^{-1}\gamma_j) = E_A(\Phi(\un_{S_{n_in_j}}) \circ \un_{S_{n_j}}^*\circ \un_{S_{n_i}})(s(\gamma_i)).$$
 But $\gamma_i \in S_{n_i}$ with $r(\gamma_i) = x \notin Y$, so
 $$E_A(\Phi(\un_{S_{n_in_j}}) \circ \un_{S_{n_j}}^*\circ \un_{S_{n_i}})(s(\gamma_i)) =
 E_A(\un_{S_{n_i}}\circ \Phi(\un_{S_{n_in_j}}) \circ \un_{S_{n_j}}^*)(r(\gamma_i)).$$
 Given $\lambda_1,\dots, \lambda_k \in \C$, we have
 $$\sum_{i,j=1}^k \lambda_i\overline{\lambda_j}F_\Phi(\gamma_{i}^{-1}\gamma_j) =
 E_A\Big(\sum_{i=1}^k (\lambda_i\un_{S_{n_i}})\circ \Phi(\un_{S_{n_i}}^* \circ \un_{S_{n_j}})\circ \sum_{j=1}^k (\lambda_j\un_{S_{n_j}})^*\Big)(x)\geq 0.$$
 \end{proof}

Obviously, if $\Phi$ is unital, $F_\Phi$ takes value $1$ almost everywhere on $X$. 

\begin{lem}\label{phi_compact} We now assume that $\Phi$ is unital, with $E_A\circ \Phi = E_A$ and $\wh{\Phi} \in \cK(\scal{M,e_A})$.
Then, for every $\varepsilon >0$, we have 
$$\nu(\set{\abs{F_\Phi} > \varepsilon} < +\infty.$$
\end{lem}

\begin{proof} Let $(S_n)$ be a partition of $G$ into Borel bisections. Given $\varepsilon >0$ we choose $\xi_1,\dots,\xi_k, \eta_1,\dots, \eta_k \in \cL^2(M,\varphi)$
such that
$$\norm{\wh{\Phi} - \sum_{i=1}^k L_{\xi_i}L_{\eta_i}^*} \leq \varepsilon/2.$$
We view $\wh{\Phi} - \sum_{i=1}^k L_{\xi_i}L_{\eta_i}^*$ as an element of $\cB(\cL^2(M,\varphi)_A)$ and we apply it to $\un_{S_n} \in \cL^2(M,\varphi)$.
Then
$$\norm{\Phi(\un_{S_n}) - \sum_{i=1}^k  \xi_i\scal{\eta_i,\un_{S_n}}_A}_{\cL^2(M)} \leq \varepsilon/2 \norm{\un_{S_n}}_{\cL^2(M)} \leq \varepsilon/2.$$

Using the Cauchy-Schwarz  inequality $\scal{\xi,\eta}_{A}^* \scal{\xi,\eta}_A \leq \norm{\xi}_{\cL^2(M)}^2\scal{\eta,\eta}_A$, we get
\begin{align*}
\norm{\scal{\un_{S_n}^*, \Phi(\un_{S_n}) - \sum_{i=1}^k  \xi_i\scal{\eta_i,\un_{S_n}}_A}_A} &\leq \norm{ \Phi(\un_{S_n}) - \sum_{i=1}^k  \xi_i\scal{\eta_i,\un_{S_n}}_A}_{\cL^2(M)}\\
&\leq \varepsilon/2.
\end{align*}

We have, for almost every $\gamma\in S_n$ and $x=s(\gamma)$,
\begin{align*}
\abs{F_\Phi(\gamma)}  &= \abs{E_A(\Phi(\un_{S_n})\circ \un_{S_n}^*)(r(\gamma))}\\
&= \abs{E_A(\un_{S_n}\circ \Phi(\un_{S_n}))(x)} =\abs{ \scal{\un_{S_n}, \Phi(\un_{S_n})}_A(x)}\\
&\leq \abs{\scal{\un_{S_n}, \Phi(\un_{S_n}) - \sum_{i=1}^k  \xi_i\scal{\eta_i,\un_{S_n}}_A}_A(x)}\\
&\quad\quad \quad\quad\quad\quad\quad+ \sum_{i=1}^k \abs{\scal{\un_{S_n},\xi_i\scal{\eta_i,\un_{S_n}}_A}_A(x)}.
\end{align*}
The first term is $\leq \varepsilon/2$ for almost every $x\in s(S_n)$. As for the second term, we have, almost everywhere,
$$\abs{\scal{\un_{S_n}, \xi_i}_A(x) \scal{\eta_i, \un_{S_n}}_A(x)} \leq \norm{\xi_i}_{\cL^2(M)} \abs{\scal{\eta_i,\un_{S_n}}_A(x)}.$$

Hence, we get
\begin{equation}\label{eq1}
\abs{F_\Phi(\gamma)} \leq \varepsilon/2 + \sum_{i=1}^k \norm{\xi_i}_{\cL^2(M)} \abs{\scal{\eta_i,\un_{S_n}}_A(s(\gamma))}
\end{equation}
for almost every $\gamma\in S_n$.

We want to estimate
$$\nu(\set{\abs{F_\Phi} > \varepsilon}) = \sum_n \nu(\set{\gamma\in S_n : \abs{F_\Phi(\gamma)} > \varepsilon}).$$
For almost every  $\gamma\in S_n$  such that $\abs{F_\Phi(\gamma)} > \varepsilon$, we see that
$$\sum_{i=1}^k \norm{\xi_i}_{\cL^2(M)} \abs{\scal{\eta_i,\un_{S_n}}_A(s(\gamma))} > \varepsilon/2.$$
Therefore
\begin{align*}
\nu(\set{\abs{F_\Phi} > \varepsilon}) &\leq \sum_n \nu\Big(\set{\gamma\in S_n : \sum_{i=1}^k \norm{\xi_i}_{\cL^2(M)} \abs{\scal{\eta_i,\un_{S_n}}_A(s(\gamma))} > \varepsilon/2}\Big)\\
&\leq \sum_n \mu\Big(\set{ x\in s(S_n) : \sum_{i=1}^k \norm{\xi_i}_{\cL^2(M)} \abs{\scal{\eta_i,\un_{S_n}}_A(x)} > \varepsilon/2}\Big).
\end{align*}
Now,
$$\sum_{i=1}^k \norm{\xi_i}_{\cL^2(M)} \abs{\scal{\eta_i,\un_{S_n}}_A(x)} \leq (\sum_{i=1}^k \norm{\xi_i}_{\cL^2(M)}^2)^{1/2} \big( \sum_{i=1}^k \abs{\scal{\eta_i,\un_{S_n}}_A(x)}^2\big)^{1/2}.$$
We set $c = \sum_{i=1}^k \norm{\xi_i}_{\cL^2(M)}^2)$ and $f_n(x) = \sum_{i=1}^k \abs{\scal{\eta_i,\un_{S_n}}_A(x)}^2$.
We have
\begin{align*}
\sum_n f_n(x) & = \sum_n \sum_{i=1}^k \abs{\scal{\eta_i,\un_{S_n}}_A(x)}^2\\
& = \sum_{i=1}^k \sum_n  \abs{\scal{\eta_i,\un_{S_n}}_A(x)}^2\\
& = \sum_{i=1}^k \scal{\eta_i, \eta_i}_A(x) \leq \sum_{i=1}^k \norm{\eta_i}_{\cL^2(M)}^2,
\end{align*}
since, by lemma \ref{calcul_scal} (or directly here),
$$
\scal{\eta_i, \eta_i}_A = \sum_k \scal{\eta_i,\un_{S_k}}_A\scal{\un_{S_k},\eta_i}
 =\sum_k \abs{\scal{\eta_i,\un_{S_k}}_A}^2.
$$

We set $d = \sum_{i=1}^k \norm{\eta_i}_{\cL^2(M)}^2$.

We have
$$\nu(\set{\abs{F_\Phi} > \varepsilon}) \leq \sum_n \mu(\set{ x\in s(S_n) : c f_n(x) > (\varepsilon/2)^2 }.$$ 
We set $\alpha = c^{-1}(\varepsilon/2)^2$. Denote by $i(x)$ the number of indices $n$ such that $f_n(x) >\alpha$.
Then $i(x) \leq N$, where $N$ is the integer part of $d/\alpha$. We denote by $\cP = \set{P_n}$ the set of subsets of $\N$
whose cardinal is $\leq N$. Then there is a partition $X = \sqcup_m B_m$ into Borel subsets such that 
$$\forall x\in B_m, \quad P_m = \set{n\in \N : f_n(x) > \alpha}.$$
We have
\begin{align*}
\nu(\set{\abs{F_\Phi} > \varepsilon}) &\leq \sum_{n,m} \mu(\set{x\in B_m\cap s(S_n) : f_n(x) > \alpha})\\
&\leq\sum_m\big(\sum_n \mu(\set{x\in B_m\cap s(S_n) : f_n(x) > \alpha})\\
&\leq\sum_m \sum_{n\in P_m} \mu(\set{x\in B_m\cap s(S_n) : f_n(x) > \alpha})\\
&\leq \sum_m N \mu(B_m) = N
\end{align*}
\end{proof}

\section{From positive type functions to completely positive maps}\label{5}
Again, we want to extend a well known result in the group case, namely that, given a positive definite function $F$ on a countable group $G$, there
is a normal completely positive map $\Phi : L(G) \to L(G)$, well defined by the formula $\Phi(u_t) = F(t) u_t$ for every $t\in G$.

We need some preliminaries. For the notion of representation used below, see definition \ref{rep}.

\begin{lem} Let $F$ be a positive definite function on $(G,\mu)$. There exists  a re\-pre\-sentation $\pi$ of $G$ on a  measurable field $\cK = \set{\cK(x)}_{x\in X}$ of
Hilbert spaces, and a measurable section $\xi : x \mapsto \xi_x\in \cK(x)$ such that 
$$F(\gamma) = \scal{\xi\circ r(\gamma), \pi(\gamma) \xi\circ s(\gamma)}$$
almost everywhere, that is $F$ is the coefficient of the representation $\pi$, associated with $\xi$.
\end{lem}

\begin{proof} This classical fact may be found in \cite{R-W}. The proof is straightforward, and similar to the classical GNS construction in the case of groups. 
Let $V(x)$ the space of
finitely supported complex-valued functions on $G^x$, endowed with the semi-definite positive hermitian form
$$\scal{f,g}_x = \sum_{\gamma_1,\gamma_2 \in G^x} \overline{f(\gamma_1)}g(\gamma_2) F(\gamma_{1}^{-1}\gamma_2).$$
We denote by $\cK(x)$ the Hilbert space obtained by separation and completion of $V(x)$, and $\pi(\gamma) : \cK(s(\gamma)) \to \cK(r(\gamma))$
is defined by $(\pi(\gamma)f)(\gamma_1) = f(\gamma^{-1}\gamma_1)$. The Borel structure on the field $\set{\cK(x)}_{x\in X}$ is provided
by the Borel functions on $G$ whose restriction to the fibres $G^x$ are finitely supported. Finally, $\mathbf{\xi}$ is the characteristic function of $X$, viewed as a Borel
section.
\end{proof}

Now we assume that $F(x) = 1$ for almost every $x\in X$, and thus $\xi$ is a unit section.
We consider the measurable field $\set{\ell^2(G_x)\otimes \cK(x)}_{x\in X}$. Note that
$$\ell^2(G_x)\otimes \cK(x) = \ell^2(G_x,\cK(x)).$$
Let $f\in \ell^2(G_x)$.  We define 
$S_x (f) \in\ell^2(G_x,\cK(x))$ by
$$S_x (f)(\gamma) = f(\gamma) \pi(\gamma)^* \xi_{r(\gamma)}$$
for $\gamma \in G_x$. Then 
$$\sum_{s(\gamma) = x} \norm{S_x (f)(\gamma)}_{\cK(x)}^2 = \norm{f}_{\ell^2(G_x)}.$$
The field $(S_x)_{x\in X}$ of operators  defines an isometry 
$$S : L^2(G,\nu)\to\int_{X}^\oplus \ell^2(G_x, \cK(x))\rd \mu(x),$$
by
$$S(f)(\gamma) = f(\gamma) \pi(\gamma)^* \xi\circ r(\gamma).$$
Note that  $\int_{X}^\oplus \ell^2(G_x, \cK(x)) \rd \mu(x)$ is a right $A$-module, by
$$(\eta a)_x = \eta_x a(x) : \gamma \in G_x \mapsto \eta(\gamma) a\circ s(\gamma).$$
Of course, $S$ commutes with the right actions of $A$.
We also observe that, as a right $A$-module,  $\cL^2(M,\varphi)\otimes_A \int_{X}^\oplus \cK(x) \rd \mu(x)$ and
 $\int_{X}^\oplus \ell^2(G_x, \cK(x))\rd\mu(x)$ are  canonically isomorphic under the map
$$\zeta\otimes \eta \mapsto \zeta \eta\circ s, \quad \forall \zeta \in \cL^2(M,\varphi), \forall\eta \in  \int_{X}^\oplus \cK(x) \rd \mu(x),$$
 where $(\zeta\eta\circ s)_x$ is the function $\gamma\in G_x \mapsto \zeta(\gamma) \eta\circ s(\gamma)$ in $\ell^2(G_x, \cK(x))$.
It follows that $M$ acts on  $\int_{X}^\oplus \ell^2(G_x, \cK(x))\rd\mu(x)$  by $m \mapsto m\otimes \Id$.
In particular, for $f \in I(G)$, we see that $L(f)\otimes \Id$, viewed as an operator on  $\int_{X}^\oplus \ell^2(G_x, \cK(x))\rd\mu(x)$ is acting as
$$((L(f)\otimes \Id)\eta)(\gamma) = \sum_{\gamma_1\gamma_2= \gamma} f(\gamma_1)\eta(\gamma_2) \in \cK(s(\gamma)).$$ 

\begin{lem}\label{construc_phi} For $f\in I(G)$, we have $S^*\big(L(f)\otimes \Id\big) S = L(F f)$.
\end{lem}
\begin{proof} A straightforward computation shows that for $\eta \in  \int_{X}^\oplus \ell^2(G_x, \cK(x))\rd\mu(x)$,
we have
$$(S^*\eta)(\gamma) = \scal{\pi(\gamma)^*\xi\circ r(\gamma), \eta(\gamma)}_{\cK(s(\gamma))}.$$ 
Moreover, given $h\in L^2(G,\nu)$, we have
\begin{align*}
\big((L(f)\otimes \Id) S h\big)(\gamma) & = \sum_{\gamma_1\gamma_2 = \gamma} f(\gamma_1)(Sh)(\gamma_2)\\
&= \sum_{\gamma_1\gamma_2 = \gamma} f(\gamma_1) h(\gamma_2) \pi(\gamma_2)^* \xi\circ r(\gamma_2).
\end{align*}
Hence,
\begin{align*}
\big(S^*(L(f)\otimes \Id)Sh\big)(\gamma) & = \scal{\pi(\gamma)^*\xi\circ r(\gamma), 
\sum_{\gamma_1\gamma_2 = \gamma} f(\gamma_1) h(\gamma_2) \pi(\gamma_2)^* \xi\circ r(\gamma_2)}\\
& = \scal{\xi\circ r(\gamma), \sum_{\gamma_1\gamma_2 = \gamma} f(\gamma_1) h(\gamma_2)  \pi(\gamma_1)\xi\circ r(\gamma_2)}\\
& = \sum_{\gamma_1\gamma_2 = \gamma} f(\gamma_1) h(\gamma_2) \scal{\xi\circ r(\gamma_1), \pi(\gamma_1)\xi\circ r(\gamma_2)}\\
&=  \sum_{\gamma_1\gamma_2 = \gamma} f(\gamma_1)F(\gamma_1) h(\gamma_2) = (L(Ff) h)(\gamma).
\end{align*} 
\end{proof}

\begin{lem}\label{F_CP} Let $F : G \to \C$ be a Borel positive type function on $G$ such that $F_{|_X} = 1$. Then there exists a 
unique normal completely positive map $\Phi$ from $M$ into $M$ such that 
$$\Phi(L(f)) = L(F f)$$
for every $f\in I(G)$. Morever, $\Phi$ is $A$-bilinear, unital and $E_A\circ \Phi = E_A$.
\end{lem}

\begin{proof} The uniqueness is a consequence of the normality of $\Phi$, combined with the density of $L(I(G))$ into $M$.
With the notation of the previous lemma, for $m\in M$ we put $\Phi(m) = S^*\big(m\otimes \Id\big) S$. Obviously, $\Phi$ satisfies the required
conditions.
\end{proof}

\begin{rem} We keep the notation of the previous lemma. A straightforward computation shows that 
$F$ is the positive definite function $F_\Phi$ constructed from $\Phi$.
\end{rem}

\begin{prop}\label{F_proper} Let $F$ be a Borel positive definite function on $G$ such that $F_{|_X} = 1$. We assume that for every $\varepsilon >0$, we have
$\nu(\set{\abs{F} > \varepsilon}) < +\infty$. Let $\Phi$ be the completely positive map defined by $F$. Then $\wh{\Phi}$ belongs to the norm closed ideal $\cI(\scal{M,e_A})$ generated by the finite projections of $\scal{M,e_A}$.
\end{prop}

\begin{proof} We observe that $T=\wh{\Phi}$ is the multiplication operator $\mathfrak{m}(F)$ by $F$. We need to show that
for every $t > 0$, the spectral projection $e_t(\abs{T})$ of $\abs{T}$ relative to $[t,+\infty[$ is finite. This projection is the multiplication
operator by $f_t = \un_{[t,+\infty[}\circ \abs{F}$.
By \eqref{formule_trace}, we have
$$\Tr_\mu(\mathfrak{m}(f_t)) =  \nu(f_t)= \nu(\set{\abs{F}> t} <+\infty.$$
\end{proof}

\section{Characterizations of the relative Haagerup property}\label{6}

\begin{thm}\label{equiv_def_property_H} The following conditions are equivalent:
\begin{itemize}
\item[(1)] $M$ has the Haagerup property relative to $A$ and $E_A$.
\item[(2)] There exists a sequence $(F_n)$ of positive definite  functions on $G$ such that
\begin{itemize}
\item[(i)] $(F_{n})_{|_X} = 1$ almost everywhere ;
\item[(ii)] for every $\varepsilon > 0$, $\nu(\set{\abs{F_n} > \varepsilon}) < +\infty$ ;
\item[(iii)] $\lim_n F_n = 1$ almost everywhere.
\end{itemize}
\end{itemize}
\end{thm}

\begin{proof} (1) $\Rightarrow$ (2). Let $(\Phi_n)$ a sequence of unital completely positive maps $M\to M$ satisfying conditions
(i), (ii), (iii) of definition \ref{def_hag_rel}. We set $F_n = F_{\Phi_n}$. By lemma \ref{phi_compact} we know that condition (ii) of (2) above is satisfied.
It remains to check (iii).  For $m\in M$, we have
$$\norm{\Phi_n(m) - m}_2^{2} = \int_X E_A((\Phi_n(m) - m)^*(\Phi_n(m) - m))(x) \rd \mu(x).$$
Let $G = \sqcup_n S_n$ be a partition of $G$ by Borel bisections. 
 There is a null subset $Y$ of $X$ such that, for every $k$ and for $\gamma\in S_k \cap r^{-1}(X\setminus Y)$ we have
 $$F_n(\gamma) - 1 = E_A(\un_{S_k}^* \circ \Phi_n(\un_{S_k}))(s(\gamma)) - E_A(\un_{S_k}^{*} \circ \un_{S_k})(s(\gamma).$$
 Thus
 \begin{align*}
 \abs{F_n(\gamma) - 1}^2 &= \abs{ E_A(\un_{S_k}^* \circ(\Phi_n(\un_{S_k}) - \un_{S_k}))(s(\gamma)}^2\\
  &\leq E_A((\Phi_n(\un_{S_k}) - \un_{S_k})^*(\Phi_n(\un_{S_k}) - \un_{S_k}))(s(\gamma).
  \end{align*}
 It follows that
   \begin{align}\label{eq_inter0}
   \int_G \abs{F_n - 1}^2 \un_{S_k} \rd \nu & \leq \int_{s(S_k)} E_A((\Phi_n(\un_{S_k}) - \un_{S_k})^*(\Phi_n(\un_{S_k}) - \un_{S_k}))(x) \rd \mu(x) \notag\\
   &\leq \norm{\Phi_n(\un_{S_k}) - \un_{S_k}}_2^{2} \to 0.
   \end{align}
   So there is a subsequence of $(\abs{F_n(\gamma) - 1})_n$ which goes to $0$ almost everywhere on $S_k$.
Using the Cantor diagonal process, we get the existence of a subsequence $(F_{n_k})_k$ of $(F_n)_n$ such that $\lim_k F_{n_k} = 1$ almost everywhere,
which is enough for our purpose.

(2) $\Rightarrow$ (1). Assume the existence of a sequence $(F_n)_n$ of positive definite  functions on $G$, satisfying the three conditions of (2).
Let $\Phi_n$ be the completely positive map defined by $F_n$. Let us show that for every $m\in M$, we have 
 $$\lim_n \norm{\Phi_n(m) - m}_2 = 0.$$
We first consider the case $m = L(f)$ with $f\in I(G)$. Then we have
$$\norm{\Phi_n(L(f)) - L(f)}_2  = \norm{L((F_n - 1)f)}_2 = \norm{(F_n - 1)f}_2 \to 0$$
by the Lebesgue dominated convergence theorem.

Let now $m\in M$. Then
$$\norm{\Phi_n(m) - m}_2 \leq \norm{\Phi_n(m-L(f))}_2 + \norm{\Phi_n(L(f))-L(f)}_2 + \norm{L(f) - m}_2.$$
We conclude by a classical approximation argument, since
$$\norm{\Phi_n(m -L(f))}_2 \leq \norm{L(f)-m}_2.$$

Together with lemmas \ref{F_CP}, \ref{F_proper} and \ref{equiv_Haag}, this proves (1).
\end{proof}
 
 This theorem justifies the following definition.
 
 \begin{defn}\label{groupoid_Haag} We say that  a countable measured groupoid  $(G,\mu)$ has the {\it Haagerup property} (or has {\it property} (H)) if 
 there exists a sequence $(F_n)$ of positive definite  functions on $G$ such that
\begin{itemize}
\item[(i)] $(F_{n})_{|_X} = 1$ almost everywhere ;
\item[(ii)] for every $\varepsilon > 0$, $\nu(\set{\abs{F_n} > \varepsilon}) < +\infty$ ;
\item[(iii)] $\lim_n F_n = 1$ almost everywhere.
\end{itemize}
\end{defn}

We observe that, by theorem \ref{equiv_def_property_H}, this notion only involves the conditional expectation $E_A$
and therefore only depends on the measure class of $\mu$. This fact does not seem to be obvious directly from
the above definition \ref{groupoid_Haag}.

\begin{defn}\label{CTN} A {\it real condi\-tionally negative definite function} on $G$ is a Borel function $\psi : G\to \R$ such that
\begin{itemize}
\item[(i)] $\psi(x) = 0$ for every $x\in G^{(0)}$;
\item[(ii)] $\psi(\gamma) = \psi(\gamma^{-1})$ for every $\gamma \in G$;
\item[(iii)] for every $x\in G^{(0)}$, every $\gamma_1,\dots,\gamma_n \in G^x$ and every real numbers $\lambda_1,\dots,\lambda_n$ with $\sum_{i=1}^n \lambda_i = 0$,
then 
$$\sum_{i,j=1}^n \lambda_i\lambda_j \psi(\gamma_{i}^{-1}\gamma_j) \leq 0.$$
\end{itemize}
\end{defn}

Such a function is non-negative.

\begin{defn}\label{CTN1} Let $(G,\mu)$ be a countable measured groupoid. A {\it real conditionally negative definite function} on $(G,\mu)$  
is a Borel function $\psi : G\to \R$ such that there exists a co-null subset $U$ of $G^{(0)}$ with the property that the restriction of $\psi$
to the inessential reduction $G_{|_U} = \set{ \gamma\in G : r(\gamma) \in U, s(\gamma)\in U}$ satisfies the conditions of the previous definition.

We say that $\psi$ is {\it proper} if for every $c>0$, we have $\nu(\set{\psi\leq c}) <+\infty$.
\end{defn}

\begin{thm} The groupoid $(G,\mu)$ has the Haagerup property if and only if there exists a  conditionally negative definite function $\psi$ on $(G,\mu)$ such that
$$\forall c > 0, \quad \nu(\set{\psi \leq c}) < + \infty.$$
\end{thm}

\begin{proof} We follow the steps of the proof given by Jolissaint \cite{Joli} for equivalence relations
and previously by Akemann-Walter \cite{A-W} for groups. Let $\psi$ be a proper conditionally negative definite function. We set $F_n = \exp(-\psi/n)$. Then $(F_n)$ is a sequence of positive definite functions
which goes to $1$ pointwise. Moreover, we have $F_n(\gamma) > c$ if and only if $\psi(\gamma) < -n\ln c$. Therefore $(G,\mu)$ has the Haagerup property.

Conversely, let $(F_n)$ be a sequence of positive definite functions on $G$ satisfying conditions (i), (ii), (iii) of theorem \ref{equiv_def_property_H} (2).
We choose sequences $(\alpha_n)$ and $(\varepsilon_n)$ of positive numbers 
such $(\alpha_n)$ is increasing with $\lim_n \alpha_n = +\infty$, $(\varepsilon_n)$ is decreasing with $\lim_n \varepsilon_n = 0$, 
and such that $\sum_n \alpha_n (\varepsilon_n)^{1/2} < +\infty$.

Let $G = \sqcup S_n$ be a partition of $G$ into Borel bisections. Taking if necessary a subsequence of $(F_n)$, we may assume,
thanks to inequality \eqref{eq_inter0},
 that for every $n$,
\begin{equation}\label{eq_inter}
\sum_{1\leq k\leq n} \int_G \abs{F_k -1}^2 \un_{S_k} \rd \nu \leq \varepsilon_{n}^2.
\end{equation}
It follows  that
\begin{align*}
\int_G \big(\Re(1-F_n\big)^2 \un_{\cup_{1\leq k \leq n} S_k}  \rd \nu
&\leq \int_G \abs{1-F_n}^2\un_{\cup_{1\leq k \leq n} S_k} \rd \nu\\
&\leq  \varepsilon_{n}^2.
\end{align*}
We set $E_n = \set{\gamma \in \cup_{1\leq k \leq n} S_k : \abs{\Re(1-F_n(\gamma)} \geq (\varepsilon_n)^{1/2}}$
and $E = \cap_l \cup_{n\geq l} E_n$. Since $\nu(E_n) \leq \varepsilon_n$ and $\sum_n \varepsilon_n < +\infty$, we see that
$\nu(E) = 0$.

Let us set $\psi = \sum_n \alpha_n \Re(1 - F_n)$ on $G\setminus E$ and $\psi = 0$ on $E$. We claim that the series converges pointwise.
Indeed, let $\gamma \in G\setminus E$. There exists $m$ such that 
$$\gamma\in (\cup_{1\leq i\leq m} S_i)\cap(\cap_{n\geq m} E_{n}^c).$$
Thus, $\abs{\Re(1-F_n(\gamma)} \leq (\varepsilon_n)^{1/2}$ for $n\geq m$, which shows our claim.

It remains to show that $\psi$ is proper. Let $c>0$, and let $\gamma\in G\setminus E$ with $\psi(\gamma) \leq c$. Then we have $\Re(1-F_n(\gamma)) \leq c/\alpha_n$
for every $n$ and therefore $\Re F_n(\gamma) \geq 1-c/\alpha_n$. Let $n$ be large enough such that $1 - c/\alpha_n \geq 1/2$. It follows that
$$\nu(\set{\psi \leq c}) \leq \nu(\set{\abs{F_n} \geq 1/2)} <+\infty.$$

\end{proof}

\begin{defn}\label{cocycle} Let $\pi$ be a representation of $(G,\mu)$ on a measurable field $\cK = \set{\cK(x)}_{x\in X}$ of Hilbert spaces. 
A {\it $\pi$-cocycle} is a Borel section $b$ of the pull-back bundle $ r : r^*\cK= \set{(\gamma, \xi) : \xi \in \cK(r(\gamma))} \to G$ such that, up to an inessential
reduction, we have, for $\gamma_1,\gamma_2 \in G$ with $s(\gamma_1) = r(\gamma_2)$,
$$b(\gamma_1\gamma_2) = b(\gamma_1) + \pi(\gamma_1) b(\gamma_2).$$

We say that $b$ is {\it proper} if for every $c >0$, we have
$\nu(\set{ \norm{b} \leq c}) < +\infty$.

\end{defn}
 
Let $b$ be a $\pi$-cocycle. It is easily seen that $\gamma\mapsto \norm{b(\gamma)}^2$ is conditionally negative definite. Moreover, every  real conditionally negative definite is of this form (see \cite[Prop. 5.21]{AD05}).

\begin{cor}\label{proper_cocycle} The groupoid $(G,\mu)$ has the Haagerup property if and only if it admits a proper $\pi$-cocycle for
some representation $\pi$.
\end{cor}

\section{Amenable countable measured groupoids have property (H)}\label{7}
Although the notion of amenable measured groupoid has been studied extensively in  \cite{AD-R}, this monograph does not
provide a characterization in term of positive definite functions. This is achieved only for topological amenability in the case of locally
compact groupoids. We now fill this gap.

Let us first recall one of the many equivalent definitions of amenability for a measured groupoid (see \cite[Prop. 3.2.14]{AD-R}).

\begin{defn}\label{def_amenable} \cite[Prop. 3.2.14 (v)]{AD-R} We say  that a countable measured groupoid $(G,\mu)$ is {\it amenable} if there exists a sequence $(\xi_n)$
of Borel functions on $G$ such that
\begin{itemize}
\item[(i)]  $\sum_{r(\gamma) = x} \abs{\xi_n(\gamma)}^2 = 1$ for almost every $x\in X$;
\item[(ii)] setting $F_n(\gamma) = \sum_{r(\gamma_1) = r(\gamma)} \overline{\xi_n(\gamma_1)}\xi_n(\gamma^{-1}\gamma_1)$, then
$\lim_n F_n = 1$  in the weak* topology on $L^\infty(G)$.
\end{itemize}
\end{defn}

This means that the restriction of $\xi_n$ to $G^x$ is a unit vector of $\ell^2(G^x)$ and that $F_n$ is the coefficient associated with $\xi_n$ of the representation
$L_H$ \footnote{We choose here the Hahn version of the regular representation (see remark \ref{Hahn}). The reason is that in this version,
the operator $\rho$ introduced in the proof of theorem \ref{carac_amen} is decomposable. } on $(\ell^2(G^x))_{x\in X}$ defined by   
$$L_H(\gamma) : \ell^2(G^{s(\gamma)}) \to \ell^2(G^{r(\gamma)}),\quad (L_H(\gamma)\xi)(\gamma_1) = \xi(\gamma^{-1}\gamma_1).$$
In order to state our equivalent characterization we introduce the following definitions. First, a Borel subset $Q$ of $G$ is said to be {\it bounded}
if there exists $c>0$ such that, for almost every $x\in X$,
$$\sharp (Q\cap G^x) \leq c \quad\hbox{and}\quad  \sharp (Q\cap G_x) \leq c.$$
We say that a Borel function $F$ on $G$ has a {\it bounded support}
if there is a bounded Borel subset $Q$ of $G$ such that $F = 0$ outside $Q$.

\begin{thm}\label{carac_amen} The groupoid $(G,\mu)$ is amenable if and only if there exists a sequence $(F_n)$ of positive definite Borel functions on $G$ such that
\begin{itemize}
\item[(i)] $(F_{n})_{|_X} = 1$ almost everywhere ;
\item[(ii)] $F_n$ has a bounded support ;
\item[(iii)] $\lim_n F_n = 1$ almost everywhere.
\end{itemize}
\end{thm}

\begin{proof} Assume  that $(G,\mu)$ is amenable. First we observe that the $\xi_n$ in definition \ref{def_amenable} may be assumed
to be non-negative (see \cite[Prop. 3.1.25, Prop. 2.2.7]{AD-R}) and so $0\leq F_n \leq 1$ almost everywhere. Let $g$ be a strictly positive $\nu$-integrable function
on $G$. We have $\lim_n \int_G g (1-F_n) \rd \nu = 0$. Hence, taking if necessary a subsequence, we may assume that $\lim_n F_n =1$
almost everywhere on $G$. It is also easy to see that, by approximation, we might have taken the $\xi_n$ to have a bounded support so that the $F_n$ too
have a bounded support.

To prove the converse, we claim that any Borel positive definite   function $F$ satisfying conditions (i), (ii), (iii) of the statement is a coefficient of the 
regular representation $L_H$ of $(G,\mu)$, that is, there exists a Borel function $\xi_F$ on $G$ such that (almost everywhere),
$$F(\gamma) = \sum_{r(\gamma_1) = r(\gamma)} \overline{\xi_F(\gamma_1)}\xi_F(\gamma^{-1}\gamma_1).$$
Since $\abs{F}$ has a bounded support and $\abs{F}\leq1$ there exists $c>0$ such that 
$$\sum_{\set{\gamma_1 : r(\gamma_1)= r(\gamma)}}\abs{F(\gamma^{-1}\gamma_1)} \leq c$$
for almost every $\gamma\in G$.

For $\xi \in L^2(G,\nu^{-1})$, we set $(\rho(F)\xi)(\gamma) = \sum_{r(\gamma_1) = r(\gamma)} F(\gamma^{-1}\gamma_1)\xi(\gamma_1)$.
We define in this way a positive operator on $L^2(G,\nu^{-1})$. Indeed, using twice the Cauchy-Schwarz inequality, we get, for $\xi,\eta\in L^2(G,\nu^{-1})$,
\begin{align*}
\scal{\abs{\eta},(\rho(\abs{F})\abs{\xi}}_{L^2(G,\nu^{-1})} &\leq \big(\int_X (\sum_{r(\gamma)=r(\gamma_1)=x}\abs{\eta(\gamma)}^2 \abs{F(\gamma^{-1}\gamma_1)})\rd \mu(x)\big)^{1/2}\\
&\quad\quad\quad\quad\quad\quad\big(\int_X (\sum_{r(\gamma)=r(\gamma_1)=x}\abs{\xi(\gamma_1)}^2 \abs{F(\gamma^{-1}\gamma_1)})\rd \mu(x)\big)^{1/2}\\
&\leq c \norm{\eta}_{L^2(G,\nu^{-1})}\norm{\xi}_{L^2(G,\nu^{-1})}.
\end{align*}
The positivity of $\rho$ is an immediate consequence of the fact that $F$ is positive definite.

It is also easy to check that $\rho$ is a decomposable operator on $L^2(G,\nu^{-1}) = \int_X^{\oplus} \ell^2(G^x) \rd \mu(x)$, and that $\rho$ commutes
with the operators $L_H(f)$. 

We set $\xi_F = \rho(F)^{1/2}\un_X$ and we write $\xi_F = \int_X^{\oplus} (\xi_{F})_x \rd \mu(x)$. A straightforward computation shows that
$$\sum_{r(\gamma) = x} \abs{\xi_F(\gamma)}^2 = F(x) = 1$$
almost everywhere. Moreover, we have
\begin{align*}
\sum_{\set{\gamma_1 : r(\gamma_1) = r(\gamma)}}\overline{\xi_F(\gamma_1)}\xi_F(\gamma^{-1}\gamma_1)
&= \scal{(\rho(F)^{1/2}\un_X)_{r(\gamma)}, L_H(\gamma)(\rho(F)^{1/2}\un_X)_{s(\gamma)}}_{\ell^2(G^{r(\gamma)})}\\
&=\scal{(\un_X)_{r(\gamma)}, L_H(\gamma)(\rho(F)\un_X)_{s(\gamma)}}_{\ell^2(G^{r(\gamma)})}\\
&= F(\gamma).
\end{align*}

This shows our claim. It follows that the existence of a sequence $(F_n)$ as in the statement of the theorem implies
the amenability of $(G,\mu)$.
\end{proof}

\begin{cor} For a countable measured groupoid, amenability implies property (H).
\end{cor}

\begin{rem} When $(G,\mu)$ is a countable measured equivalence relation, theorem \ref{carac_amen} is to compare with
 the deep result of Connes-Feldman-Weiss \cite{CFW}
saying that, whenever amenable, this relation is hyperfinite : there is an increasing sequence $(\cR_n)$ of finite subequivalence relations
such that $\cup_n \cR = G$ almost everywhere. Then the characteristic function $F_n$  of $\cR_n$ is positive definite (see lemma \ref{ext_def_pos}),
and $\lim _n F_n = 1$ almost everywhere. It is easily seen that we may take the $\cR_n$ to be bounded.
\end{rem}

\section{Treeable countable measured groupoids have property (H)}\label{8}
The notion of treeable countable measured equivalence relation has been introduced by Adams in \cite{Ad}. Its obvious
extension to the case of countable measured groupoids is exposed in \cite{AD05}. We recall here the main definitions.
Let $Q$ be a Borel subset of a countable Borel groupoid  $G$. We set $Q^0 = X$ and for $n\geq 1$, we set
$$Q^n = \set{\gamma\in G : \exists \gamma_1, \dots, \gamma_n \in Q, \gamma = \gamma_1\cdots\gamma_n}.$$
\begin{defn} A {\it graphing} of $G$ is a Borel subset $Q$ of $G$ such that $Q = Q^{-1}$, $Q\cap X = \emptyset$ and $\cup_n Q^n = G$.
\end{defn}

A graphing defines a structure of $G$-bundles of graphs (definition \ref{graphs}), where $G$ is the set of vertices and
$$\E = \set{(\gamma_1,\gamma_2)\in G\,_r\!*_r G :  \gamma_{1}^{-1}\gamma_2 \in Q}$$
is the set of edges. Thus, a graphing is a Borel way of defining a structure of graph on each fibre $G^x$. These graphs are
connected since $\cup_n Q^n = G$.

Note that $G^x$ is endowed with the length metric $d_x$ defined by
$$d_x( \gamma_1,\gamma_2) = \min\set{n\in \N : \gamma_1^{-1}\gamma_2 \in Q^n}.$$
The map $(\gamma_1,\gamma_2) \in G\,_r\!*_r G= \set{(\gamma_1,\gamma_2) : r(\gamma_1) = r(\gamma_2)} \mapsto d_{r(\gamma_1)}(\gamma_1,\gamma_2)$ is Borel.
It follows that a graphing defines on $G$  a structure of Borel $G$-bundle of metric spaces (see definition \ref{G_metric}). 

The graphing  $Q$ is called a {\it treeing} if it gives to $G$ the structure of a $G$-bundle of trees (see definition \ref{G_tree}).

\begin{defn}\label{def_treeable}  A countable Borel groupoid $G$ is said to be {\it treeable} if there is a graphing which gives to $r: G\to X$ a structure
of  $G$-bundle of trees.

A countable measured groupoid $(G,\mu)$ is said to be {\it treeable} if there exists an inessential reduction $G_{|_U}$ which is a treeable Borel
groupoid in the above sense.

Equipped with such a structure, $(G,\mu)$ is said to be a {\it treed measured groupoid}.
\end{defn}

Consider the case where $G$ is a countable group and $Q$ is a symmetric set of generators. The corresponding graph structure on $G$
is the Cayley graph defined by $Q$. If $Q = S \cup S^{-1}$ with $S\cap S^{-1} = \emptyset$, then $Q$ is a treeing if and only if $S$
is a free subset of generators of $G$ (and thus $G$ is a free group).

We warn the reader that our notion of treeing is slightly different from the notion used by Levitt and Gaboriau for equivalence relations.
The nuance is analyzed in \cite{Gab_cout} and \cite[\S 17]{KM}. In particular the treeings in Levitt sense are provided with a natural orientation
(in the sense of  definition \ref{G_tree}), which
is not the case for us. However, every treed equivalence relation is orientable. For general
groupoids, this is almost true, as shown in the following lemma.

 \begin{lem}\label{lem_orien_tre} Let $(G,\mu)$ be a countable measured groupoid and let $Q$ be a treeing of $G$. We assume that for almost every $x\in X$, 
 the set $Q \cap G(x)$ does not contain elements of period $ 2$. Then  $Q$ is orientable.
 \end{lem} 
 \begin{proof}
Let $f: Q \to [0,1]$ be a Borel injection. We set 
$$Q_+ = \set{\gamma \in Q : f(\gamma) < f(\gamma^{-1})}.$$
We have $Q = Q_+ \cup Q_+^{-1}$ and $\emptyset = Q_+ \cap Q_+^{-1}$.
 Then 
 $$\E_+ = \set{(\gamma_1,\gamma_2)\in G^x\times G^x : \gamma_{1}^{-1}\gamma_2 \in Q_+}$$
 provides an orientation of the $G$-bundle of trees.

\end{proof}

As made precise  in \cite[Prop. 3.9]{Alv}, treeable groupoids are the analogue of free groups and therefore the following theorem is no surprise.

\begin{thm}[Ueda]\label{Ueda} Let $(G,\mu)$ be a countable measured groupoid which is treeable. Then $(G,\mu)$ has
the Haagerup property.
\end{thm}

Let $Q$ be a treeing of $(G,\mu)$ and let $(d_x)_{x\in X} $ be the associated field of metrics $d_x$ on $G^x$.
We set $\psi(\gamma) = d_{r(\gamma)}(r(\gamma), \gamma)$. It is a real conditionally definite negative function on $G$.
Indeed, given $\gamma_1,\dots,\gamma_n \in G^x$ and $\lambda_1,\dots,\lambda_n \in \R$ such that $\sum_{i=1}^n \lambda_i = 0$, we have
\begin{align*}
\sum_{i,j = 1}^n \lambda_i \lambda_j \psi(\gamma_i^{-1}\gamma_j) &= \sum_{i,j = 1}^n \lambda_i \lambda_j d_{s(\gamma_i)}(s(\gamma_i),\gamma_i^{-1}\gamma_j)\\
& = \sum_{i,j = 1}^n \lambda_i \lambda_j d_x(\gamma_i,\gamma_j) \leq 0,
\end{align*}
since the  length metric on a tree is conditionally definite negative (see \cite[page 69]{HV} for instance).

We begin by proving theorem \ref{Ueda} in the case where $Q$ is bounded.

\begin{lem}\label{tree_inter} Assume that $Q$ is bounded. 
Then, for every $c>0$ we have $\nu(\set{\psi  \leq c}) < +\infty$.
\end{lem}

\begin{proof} We have
$$\nu(\set{\psi  \leq c}) = \int_X \sharp\set{\gamma : s(\gamma) = x, d_x(x,\gamma^{-1})\leq c} \rd\mu(x).$$
Let $k$ be such that $\sharp Q^x \leq k$ for almost every $x\in X$. The cardinal of the ball in $G^x$ of center $x$ and radius $c$
is smaller than $k^c$. It follows that $\nu(\set{\psi  \leq c}) \leq k^c$.
\end{proof}

In view of the proof in the general case, we make a preliminary observation.
 Whenever $Q$ is bounded, $G$ is the union of the increasing sequence $(\set{\psi \leq k})_{k\in \N}$ of Borel subsets with $\nu(\set{\psi \leq k})<+\infty$.
Then, setting $F_n = \exp(-\psi/n)$, we have  $\lim_n F_n = 1$ uniformly on each subset $\set{\psi \leq k}$. Indeed, if $\psi(\gamma) \leq k$, we get
$$0\leq 1- F_n(\gamma) \leq \sum_{j\geq 1} \frac{1}{n^j} \frac{\psi(\gamma)^j}{j!}
\leq \frac{k}{n} \exp{(k/n)}.$$

\begin{proof}[Proof of theorem \ref{Ueda}] The treeing $Q$ is no longer supposed to be bounded. Let $G = \sqcup S_k$ be a partition of $G$
into Borel bisections. For every integer $n$ we set
$$Q_{n}' = \cup_{k\leq n} (Q\cap S_k)\quad\hbox{and}\quad Q_n = Q_{n}' \cup (Q_{n}')^{-1}.$$
Note that $(Q_n)$ is an increasing sequence of Borel symmetric and bounded subsets of $Q$ with $\cup_n Q_n = Q$.
Let $G_n$ be the subgroupoid of $G$ generated by $Q_n$, that is $G_n = \cup_{k\geq 0} Q_{n}^k$, where we put $Q_{n}^0 = X$.

We observe that $Q_n$ is a treeing for $G_n$. Denote by $\psi_n$ the associated conditionally definite negative function on $G_n$.
Since $Q_{n-1} \subset Q_n$, we have 
$$(\psi_{n})_{|_{G_{n-1}}} \leq \psi_{n-1}.$$
Given two integers $k$ and $N$, we set
$$A_{k,N} = \set{\gamma\in G_k : \psi_k(\gamma) \leq N}.$$
Then, obviously we have
$$A_{k,N}\subset A_{k+1,N}\quad\hbox{and}\quad A_{k,N}\subset A_{k,N+1}.$$
In particular, $(A_{k,k})_k$is an increasing sequence of Borel subsets of $G$ with $\cup_k A_{k,k} = G$.

We fix $k$. We set $F_{k,n}(\gamma) = \exp(-\psi_k(\gamma)/n)$ if $\gamma\in G_k$ and $F_{k,n}(\gamma) = 0$ if $\gamma\notin G_k$.
By lemma \ref{ext_def_pos} to follow, $F_{k,n}$ is positive definite on $G$. Since $Q_k$ is bounded, lemma \ref{tree_inter} implies that for every
$\varepsilon >0$, and for every $n$, we have $\nu(\set{F_{k,n} \geq \varepsilon}) <+\infty$. Moreover,  $\lim_n F_{k,n} = 1$
uniformly on each $A_{k,N}$, $N\geq 1$, as previously noticed.

We  choose, step by step, a strictly increasing sequence $(n_i)_{i\geq 1}$ of integers such that for every $k$,
$$\sup_{\gamma\in A_{k,k}} 1- F_{k,n_k}(\gamma) \leq 1/k.$$
Then the sequence $(F_{k,n_k})_k$ of positive definite functions satisfies the required conditions showing that $(G,\mu)$ has property (H).
\end{proof}

\begin{lem}\label{ext_def_pos} Let $H$ be a subgroupoid of a groupoid $G$ with $G^{(0)} = H^{(0)}$. Let $F$ be a positive definite function
on $H$ and extend $F$ to $G$ by setting $F(\gamma) = 0$ if $\gamma\notin H$. Then $F$ is positive definite on $G$.
\end{lem}

\begin{proof} Let $\gamma_1,\dots,\gamma_n\in G^x$ and let $\lambda_1,\dots,\lambda_n \in \C$. We want to show that
\begin{equation*}\label{to_show}
\sum_{i,j=1}^n \lambda_i \overline{\lambda_j} F(\gamma_i^{-1}\gamma_j) \geq 0.
\end{equation*}
We assume that this inequality holds for every number $k<n$ of indices.
For $k=n$, this inequality  is obvious if for every $i\not=j$ we have $\gamma_i^{-1}\gamma_j \notin H$. Otherwise, up to a permutation of indices, 
 we  take $j = 1$ and we assume that $2,\dots, l$ are the indices $i$ such that $\gamma_i^{-1}\gamma_1 \in H$.
Then, if $1\leq i,j \leq l$ we have $\gamma_{i}^{-1} \gamma_j = (\gamma_{i}^{-1} \gamma_1)(\gamma_{1}^{-1}\gamma_j) \in H$
and for $i\leq l < j$ we have $\gamma_{i}^{-1} \gamma_j \notin H$. It follows that
$$\sum_{i,j=1}^n \lambda_i \overline{\lambda_j} F(\gamma_i^{-1}\gamma_j) = \sum_{i,j=1}^l \lambda_i \overline{\lambda_j} F(\gamma_i^{-1}\gamma_j)
+ \sum_{i,j> l} \lambda_i \overline{\lambda_j} F(\gamma_i^{-1}\gamma_j),$$
where the first term of the right hand side is $\geq 0$. As for the second term, it is also $\geq 0$
by the induction assumption.
\end{proof}

\section{Properties (T) and (H) are not compatible}\label{9}

Property (T) for group actions and equivalence relations has been introduced  by Zimmer in \cite{Zim81}. Its extension to measured groupoids
is immediate. We say that $(G,\mu)$ has {\it property (T)} if whenever a representation of $(G,\mu)$ almost has unit invariant section, it actually
has a unit invariant section (see \cite[def. 4.2, def. 4.3]{AD05} for details). We have proved in \cite[thm. 5.22]{AD05} the following characterization of property (T).

\begin{thm} Let $(G,\mu)$ be an ergodic countable measured groupoid. The follo\-wing conditions are equivalent:
\begin{itemize}
\item[(i)] $(G,\mu)$ has property (T);
\item[(ii)] for every real conditionally definite negative function $\psi$ on $G$, there exists a Borel subset $E$ of $X$, with $\mu(E) >0$,
such that the restriction of $\psi$ to $G_{|_E}$ is bounded.
\end{itemize}
\end{thm}

\begin{thm}\label{H-T} Let $(G,\mu)$ be an ergodic countable measured groupoid. We assume that $(G^{(0)},\mu)$ is a diffuse standard probability space.
Then $(G,\mu)$ cannot have simultaneously properties (T) and (H).
\end{thm}

\begin{proof} Assume that $(G,\mu)$ has both properties (H) and (T). There exists  a Borel conditionally definite negative function  $\psi$ such
that for every $c > 0$, we have $\nu(\set{\psi \leq c}) < +\infty$. Moreover, there exists a Borel subset $E$ of $X$, with $\mu(E) >0$,
such that the restriction of $\psi$ to $G_{|_E}$ is bounded, say by $c$. Then, we have
$$\int_E \sharp\set{\gamma : s(\gamma) = x, r(\gamma)\in E} \rd\mu(x) <+\infty.$$
Therefore, for almost every $x\in E$, we have $\sharp\set{\gamma : s(\gamma) = x, r(\gamma)\in E} < +\infty$. Replacing if  necessary $E$
by a smaller subset we may assume the existence  of $N>0$ such that for every $x\in E$,
$$\sharp\set{\gamma : s(\gamma) = x, r(\gamma)\in E}\leq N.$$

Since $(G_{|_E}, \mu_{|_E})$ is ergodic, we may assume that all the fibres of this groupoid have the same finite cardinal. Therefore,
this groupoid is proper and so the quotient Borel space $E/(G_{|_E})$ is countably separated (see \cite[lemma 2.1.3]{AD-R}). A classical
argument (see \cite[prop. 2.1.10]{Zim_book}) shows that $\mu_{|_E}$ is supported by an equivalence class, that is by a finite subset of $E$.
But this contradicts the fact that the measure is diffuse.
\end{proof}

\begin{prop} Let $(G,\mu)$ be a countable ergodic measured groupoid such that $(\cR_G,\mu)$ has property (H) (for instance is treeable). We assume that $(X,\mu)$ is a diffuse
probability measure space. Then $(G,\mu)$ has not property (T).
\end{prop}

\begin{proof} If $(G,\mu)$ had property (T) then $(\cR_G,\mu)$ would have the same property by \cite[thm. 5.18]{AD05}. But this is impossible by theorem
\ref{H-T}.
\end{proof}

This allows to retrieve results of Jolissaint \cite[prop. 3.2]{Joli} and Adams-Spatzier \cite[thm. 1.8]{AS}.

\begin{cor}\label{T-non-H} Let $\Gamma \actson (X,\mu)$ be an ergodic  probability measure preserving action of a countable group $\Gamma$ having property (T).
Then $(\cR_\Gamma,\mu)$ has not property (H) and in particular is not treeable.
\end{cor}

\begin{proof} Indeed under the assumptions of the corollary, the semi-direct product groupoid $(X\rtimes \Gamma, \mu)$ has property (T)
by \cite[prop. 2.4]{Zim81}, and we apply the previous theorem.
\end{proof}

\section{Further remarks}\label{10}
As already mentioned, the groupoid $(G,\mu)$ is an extension of the associated equivalence relation $(\cR_G,\mu)$ by the {\it stabilizer groupoid}, that
is the bundle $(G(x))_{x\in X}$
of isotropy groups. We review how the different properties of $(G,\mu)$ studied in the previous sections are related to the corresponding properties 
for $(\cR_G,\mu)$ and the isotropy groups.

\subsection{Amenability}
 
It is known (see \cite{GS,AEG} for group actions and \cite[cor. 5.3.33]{AD-R} in the general case) that $(G,\mu)$ is amenable if
and only if $(\cR_G,\mu)$ is amenable and almost every isotropy group $G(x)$ is amenable.

\subsection{Property (H)} Let $\Gamma \actson X$ be an action of countable group, which preserves the class of a probability measure $\mu$.
Obviously, the semi-direct product groupoid $(G = X \rtimes\Gamma, \mu)$ has property (H) whenever $\Gamma$ has property (H). Indeed, let
$(f_n)$ be a sequence of positive definite functions vanishing at infinity on $\Gamma$, going to $1$ pointwise. If we set $F_n(x,g) = f_n(g)$, we get
a sequence of positive definite functions on $G$ satisfying the properties of definition \ref{groupoid_Haag}. The converse is true for a free
probability measure preserving action (see \cite[Prop. 3.5]{Joli02} or \cite[Prop. 3.1]{Popa06}).

Property (H) for $(G,\mu)$ does not imply that $(\cR_G,\mu)$ has property (H).
Indeed consider a free probability measure preserving action of a countable group $\Gamma$ on $(X,\mu)$,
where $\Gamma$ has not property (H). We write $\Gamma$ as a quotient of a free group $F$ and we let $F$ act on $X$ through $\Gamma$.
Then $(G=X\rtimes F, \mu)$ has property (H) although  the associated equivalence relation $(\cR_G,\mu)  = (\cR_{X\rtimes\Gamma},\mu)$
does not share this property. Note that $(X \rtimes F,\mu)$ is even a treeable groupoid.

\begin{prop}\label{isot_haag} Let $(G,\mu)$ be a countable measured groupoid with property (H). Then almost every isotropy group $G(x)$  has property (H).
\end{prop}

\begin{proof} Let $\psi$ be a Borel conditionally negative definite function such that for every $c>0$,
$$\int_X \sharp(G_x \cap \set{\psi \leq c}) \rd\mu(x) = \nu(\set{\psi \leq c}) <+\infty.$$
Then there is a conull subset $U$ of $X$ such that for every $x\in U$ and every integer $n$ we have
$$\sharp(G_x \cap \set{\psi \leq n}) < +\infty.$$
By considering the restriction of $\psi$ to the group $G(x) \subset G_x$ we deduce that $G(x)$ has property (H) whenever $x\in U$.
\end{proof}
\subsection{Treeability} We have  given above an example of a treeable groupoid such that the associated equivalence relation has not property (H).
A fortiori, it is not treeable.

\begin{prop} Let $(G,\mu)$ be an orientable treeable countable measured groupoid. Then almost every isotropy group $G(x)$ is a free group.
\end{prop}

\begin{proof} By assumption, passing to an inessential reduction, there is a structure of $G$-bundle of trees for $r : G \to X$. Each $G(x)$ acts freely and without inversion on the tree $G^x$ and so is a free group by \cite[thm. 4, p. 41]{Serre}.
\end{proof} 

\begin{rem} Let $(G,\mu)$ be an ergodic amenable countable measured groupoid.  Then $(\cR_G,\mu)$
is an ergodic amenable equivalence relation and therefore isomorphic to an equivalence relation defined by a free $\Z$-action. It is therefore treeable.
Assume in addition that  $(G,\mu)$ has an orientable treeing.
Then almost every isotropy group is amenable and free, so isomorphic to $\Z$ or trivial. 

Conversely, if $(\cR_G,\mu)$ is amenable and if almost every isotropy group is isomorphic to $\Z$ or trivial, then $(G,\mu)$ is amenable.
Is is treeable ?
\end{rem}

\subsection{Property (T)}
If $(G,\mu)$ is an ergodic countable groupoid having property (T), then $(\cR_G,\mu)$ too has property (T) (see \cite[thm. 5.18]{AD05}).
On the other hand, the isotropy groups have not always   property (T) (almost everywhere). Indeed, consider the obvious action of $SL_3(\Z)$
on $\T^3$. It is essentially free and preserves the Haar measure $\mu$. We let the group $\Z^3 \rtimes SL_3(\Z)$ act on
$\T^2$ through $SL_3(\Z)$. Since $\Gamma =\Z^3 \rtimes SL_3(\Z)$ has property (T) and since the action is probability measure preserving,
the groupoid $(\T^3 \rtimes \Gamma,\mu)$ has property (T) (by \cite[Prop. 2.4]{Zim81}). However, the isotropy groups for this action
are $\Z^3$ almost everywhere.

It is known that any group extension of property (T) groups keeps this property. It is likely that a countable measured groupoid
whose associated equivalence relation and (a.e.) isotropy groups have property (T) has this property.

\section{Appendix : $G$-bundles} \label{11} For the reader's convenience, we recall some basic definitions. For more details we refer
to \cite{AD05}.

A {\it bundle over a Borel space } $Y$ is a Borel space $Z$ equipped with a Borel surjection $p : Z\to Y$. The Borel spaces $Z_y = p^{-1}(y)$
are the fibres of the bundle. Given another bundle $p' : Z' \to Y$, the {\it fibred product} of $Z$ and $Z'$ over $Y$ is the bundle
$$Z*Z' = \set{(z,z') \in Z\times Z' : p(z) = p'(z)}.$$
In case of ambiguity, we shall write $Z\,_{p}\!*_{p'} Z'$ instead of $Z*Z'$.

 A {\it section} of the bundle $p: Z\to Y$ is a Borel map $\xi : Y \to Z$ with
$p\circ\xi(y) = y$ for all $y$.

\begin{defn}\label{G_bundle} Let $G$ be a Borel groupoid. A (left) {\it $G$-bundle}, or (left) {\it $G$-space}, is a bundle $p: Z \to X= G^{(0)}$
endowed with a Borel map $(\gamma,z) \mapsto \gamma z$ from $G\,_s *_p Z$ to $Z$ such that, whenever the operations make sense,
\begin{itemize}
\item[(i)] $p(\gamma z) = r(\gamma)$, $p(z) z = z$;
\item[(ii)] $(\gamma_1 \gamma_2) z = \gamma_1(\gamma_2 z)$ whenever the operations make sense.
\end{itemize}
 An {\it invariant section} is a  section $\xi$ such that 
$\gamma\xi(x) = \xi(\gamma x)$ for $(\gamma,x) \in G * X$.
\end{defn}

\begin{defn}\label{graphs} Let $G$ be a Borel groupoid. A (left) {\it $G$-bundle of graphs} is a pair $(\V,\E)$ of (left) $G$-bundles
where $\E$ is a subset of $\V *\V$ which does not intersect the diagonal, the $G$-action on $\E$ being given by
$$\gamma(v_1,v_2) = (\gamma v_1,\gamma v_2).$$
In particular, for every $x\in X$, the fibre $\V_x$ is a graph with $\V_x$ as set of vertices and $\E \cap( \V_x \times \V_x)$ as set of edges.

An {\it oriented $G$-bundle of graphs} is a $G$-bundle $(\V,\E)$ of graphs together
with a Borel subset $\E^+$ of $\E$, with $\E^+ \cap \overline{\E}^+ = \emptyset$,
$\E = \E^+ \cup \overline{\E}^+$ and $G \E^+ = \E^+$, where $\overline{\E}^+ :=
\{(z',z) : (z,z')\in \E^+\}$ denotes the set of opposite edges of $\E^+$. The subset $\E^+$
is called an {\it orientation} of the $G$-bundle of graphs.
One also says that $G$ {\it acts on the bundle of graphs without inversion}.

We say that a $G$-bundle of graphs is {\it orientable} if it may be given an orientation.
\end{defn}

\begin{defn}\label{G_tree} A {\it $G$-bundle of trees} is a  $G$-bundle $(\V,\E)$ of graphs, such that 
each fibre $\V_x$ is a countable tree.

\end{defn}

\begin{defn}\label{metri} Let $Y$ be a Borel space. A {\it bundle of metric spaces
over} $Y$ is a bundle $p: Z\to Y$ equipped with a family 
$(d_y)_{y\in Y}$  of
metrics $d_y$ on $Z_y$ satisfying the two following conditions:
\begin{itemize}
\item[(1)] the map $(z,z') \mapsto d_{p(z)}(z,z')$ defined on $Z*Z$ is Borel.
\item[(2)] there exists a sequence $(\xi_n)$ of sections such that for every
$y\in Y$, the set $\{\xi_n(y) : n\in \N\}$ is dense in $Z_y$.
\end{itemize}
\end{defn}

As a consequence of the above condition $(1)$, note that for every section $\xi$
of $Z$, the map $(y,v) \mapsto d_y(v,\xi(y))$ is Borel on $Y*Z$.

\begin{defn}\label{G_metric} Let $G$ be a Borel groupoid. A {\it $G$-bundle of metric spaces}  is a $G$-space $p: Z \to X$ which is also a bundle
 of metric spaces such that $G$ acts fibrewise by isometries. 
\end{defn}

Let now $\cH = \{\cH(y)\}_{y \in Y}$ be a family of Hilbert spaces  indexed by a
Borel set
$Y$ and denote by $p$ the projection from $Y*\cH = \{(y,v) :v \in \cH(y)\}$ 
to $Y$.

\begin{defn}[\cite{Ram}, p. 264] A {\it Hilbert bundle} (or {\it Borel field of Hilbert spaces}) on a Borel space $Y$ is a
space
$Y*\cH$ as above, endowed with a Borel structure such that
\begin{itemize}
\item[(1)] a subset $E$ of $Y$ is Borel if and only if $p^{-1}(E)$ is Borel;
\item[(2)] there exists a {\it fondamental sequence $(\xi_n)$ of sections} satisfying
the following conditions :
\begin{itemize}
\item[(a)] for every $n$, the map $(y,v) \rightarrow \langle \xi_n(y), v\rangle$ is
Borel on $Y*\cH$;
\item[(b)] for every $m,n$, the map $y\rightarrow \langle \xi_m(y), \xi_n(y)\rangle$ is Borel;
\item[(c)] for every $y\in Y$, the set $\{\xi_n(y) : n\in \N\}$ is dense in $\cH(y)$.
\end{itemize}
\end{itemize}
If $\mu$ is a probability measure on $Y$, a {\it measurable field of Hilbert space} on $(Y,\mu)$
is a Hilbert bundle on some conull subset of $Y$.
\end{defn}

Given a Hilbert bundle $\cH$ over the set $X$ of units of a Borel groupoid $G$, we denote by
$\hbox{Iso}(X*\cH)$ the groupoid formed by the triples $(x,V,y)$, where
$x,y\in X$ and $V$ is a Hilbert isomorphism from $\cH(y)$ onto $\cH(x)$,
the composition law being defined by $(x,V,y)(y,W,z) = (x, V\circ W, z)$. We endow 
$\hbox{Iso}(X*\cH)$ with the weakest Borel structure such that
$(x,V,y) \mapsto \scal{V\xi_n(y), \xi_m(x)}$ is Borel for every $n,m$, where
$(\xi_n)$ is a fondamental sequence.

\begin{defn}\label{rep} A {\it representation of a Borel groupoid} $G$ is a pair
$(X*\cH, \pi)$ where $X*\cH$ is a Hilbert bundle over $X =G^{(0)}$,
and $\pi : G \rightarrow \hbox{Iso}(X*\cH)$ is a Borel homomorphism that 
preserves the unit space. For $\gamma\in G$, we have $\pi(\gamma)
= (r(\gamma), \hat{\pi}(\gamma), s(\gamma)$, where $\hat{\pi}(\gamma)$ is an isometry
from
$\cH (s(\gamma))$ onto  $\cH(r(\gamma))$. To lighten the notations, we shall
identify $\pi(\gamma)$ and $ \hat{\pi}(\gamma)$. In particular we have
\begin{equation}\label{repre}
\forall (\gamma_1, \gamma_2) \in G^{(2)},\, \pi(\gamma_1\gamma_2) = \pi(\gamma_1)\pi(\gamma_2);\quad
\forall \gamma\in G,\, \pi(\gamma^{-1}) = \pi(\gamma)^{-1}.
\end{equation} 
\end{defn}

Now let   $(G,\mu)$ be a  countable measured groupoid. As usual, in this setting it is enough to consider the above notions of $G$-bundles,
up to an inessential reduction. For instance a representation of $(G,\mu)$ is a representation $\pi$ of some inessential reduction $G|_U$.
\bibliographystyle{plain}

\end{document}